\documentclass[nospthms]{svjour3}
\smartqed
\usepackage{amssymb,amsfonts,tikz}
\usepackage{amsthm}
\usepackage{amsmath}
\usepackage{calc,color}
\usepackage{graphics}
\usepackage{pdfsync}
\usepackage{epstopdf}
\usepackage{enumerate}
\usepackage{bm,url}
\usepackage{cases}
\usepackage{graphics}
\usepackage{pdfsync}
\usepackage{epstopdf}
\usetikzlibrary{hobby}
\usepackage{enumerate}
\usepackage{fullpage}
\usepackage{graphicx}
\usepackage{empheq}
\usepackage{relsize}
\usepackage{esint}
\usepackage{hyperref}
\usepackage{caption}
\usepackage{subcaption}

\newcommand{\oG}{\mathcal{G}}
\newcommand{\oL}{\mathcal{L}}
\newcommand{\Int}{\displaystyle\int}

\newcommand{\re}{\mathbb{R}}
\newcommand{\na}{\mathbb{N}}
\newcommand{\ds}{\displaystyle}
\newcommand{\mc}{\mathcal}
\newcommand{\dom}{\Omega}
\newcommand{\bnd}{\Gamma}
\newcommand{\nd}{\quad\text{ and }\quad}

\newcommand{\wt}[1]{\widetilde{#1}}

\newcommand{\vep}{\varepsilon}

\def\sym{\operatorname{sym}}
\def\asym{\operatorname{asym}}

\newcommand{\bs}{\setminus}

\newcommand{\ren}{\re^n}

\newcommand{\smfrac}[2]{{\textstyle{\frac{#1}{#2}}}}

\newcommand{\lp}{\left(}
\newcommand{\rp}{\right)}
\newcommand{\ls}{\left[}
\newcommand{\rs}{\right]}
\newcommand{\lc}{\left\{}

\newcommand{\lt}{\left.}
\newcommand{\rt}{\right.}

\newtheorem{theorem}{Theorem}[section]
\newtheorem{corollary}[theorem]{Corollary}
\newtheorem{lem}[theorem]{Lemma}

\newtheorem{ex}[theorem]{Example}
\theoremstyle{definition}

\theoremstyle{remark}
\newtheorem{rem}[theorem]{Remark}
\numberwithin{equation}{section}

\title{
Sensitivity analysis for solutions to heterogeneous nonlocal systems
\thanks{N.B, M.F, P.R. were supported by the award NSF -- DMS 1716790}}

\subtitle{Theoretical and numerical studies}

\date{July 2020}
\author{Nicole E. Buczkowski \and Mikil D. Foss \and Michael L. Parks \and Petronela Radu}
\institute{Nicole Buczkowski
            \at
             University of Nebraska-Lincoln\\        \email{nbuczkowski@huskers.unl.edu}
         \and Mikil D. Foss \at
              University of Nebraska-Lincoln \\
              \email{mfoss@unl.edu}
         \and Petronela Radu
            \at
              University of Nebraska-Lincoln \\
              \email{pradu@unl.edu}
         \and Michael L. Parks
            \at
              Sandia National Laboratories \\
              \email{mlparks@sandia.gov}
}

\date{Received: date / Accepted: date}

\begin{document}

\maketitle

\begin{abstract} The paper presents  a collection of results on continuous dependence for solutions to nonlocal problems under perturbations of data and system parameters. The integral operators appearing in the systems capture interactions via heterogeneous kernels that exhibit different types of weak singularities, space dependence, even regions of zero-interaction. The stability results showcase explicit bounds involving the measure of the domain and of the interaction collar size, nonlocal Poincar\'e constant, and other parameters.  In the nonlinear setting the bounds quantify in different $L^p$ norms the sensitivity of solutions under different nonlinearity profiles. The results are validated by numerical simulations showcasing discontinuous solutions,  varying horizons of interactions, and symmetric and heterogeneous kernels.
\end{abstract}

\section{Introduction}

The third condition of Hadamard well-posedness, continuous dependence on data, is important for several reasons in mathematical models. From an application standpoint, whenever data is based on physical observations it is expected to have some associated measurement errors, due either to inaccessibility or prohibitive costs. A system where solutions depend continuously on the given data will ensure that the values of the approximated solution can be found in a close range near the exact solution. In fact, it is desirable to obtain {\it explicit bounds} that quantify the effect that noise or small perturbations in data or parameters of the system can have on solutions. In numerical simulations, a decrease in the mesh spacing is anticipated to lead to better approximations of the exact solution at a prescribed rate (i.e., numerical convergence). However, numerical convergence theorems in general apply only to well-posed mathematical problems, meaning that continuous dependence of the solution on the data is a necessary (although not sufficient) condition for numerical convergence. 

We will conduct these stability studies in the nonlocal framework which has generated interest due to its capability to handle discontinuities, both in the input functions themselves, as well as in the domains where the equations are posed. The operators have an integral form which collect information in a neighborhood of a point through a kernel of interaction. Nonlocal interactions in a variety of applications have been expressed through different operators (single, or double convolution-type operators), for which results connecting the local and nonlocal frameworks have been established \cite{alali2015generalized,andreu2010nonlocal,du2013nonlocal,gal2017strong,radu2017nonlocal,radu2019doubly,Silling,silling2010linearized}. 

The results here are formulated on a bounded domain $\Omega$, with the linear systems involving a nonlocal Laplacian operator:
\begin{equation}\label{lap}
\mc{L}_{\mu} u(x)=\Int_{\re^n} (u(y)-u(x))\mu(x,y)dy, \ x\in \Omega,
\end{equation}
where the kernel $\mu(x,y)$ that records the interaction between points $x$ and $y$ is chosen to be nonnegative and integrable (in order to allow ``rough" inputs). In peridynamics, the kernel is usually chosen to be symmetric, though for the purposes of this paper, we allow the kernel to be heterogeneous. The kernel allows an added degree of flexibility in applications, so for different profiles we obtain different dynamics. 
At a theoretical level, the selection of a weakly singular function (vs. a highly singular one) allows discontinuous solutions, while the availability of mathematical tools (such as compactness theorems) is highly reduced.
While a growing literature dedicated to nonlocal Laplacians of this form establishes growing connections between classical and elliptic-type properties \cite{foss2019bridging}, investigations regarding continuous dependence have been more scarce. Several studies focus on nonlocal models which include fractional operators (where kernels have nonintegrable singularity); we summarize some of the results below.  Sensitivity on system parameters (such as operator coefficients) has been explored for systems in \cite{burkovska2020affine} and \cite{tuan2020continuity}. Continuous dependence on initial conditions for nonlinear fractional convection-diffusion was studied in  \cite{alibaud2012continuous} and for nonlinear nonlocal diffusion in \cite{bogoya2008nonlocal}. The latter also includes some continuous dependence on boundary conditions. The paper \cite{coclite2018wellposedness} includes results on stability for boundary or initial data, where the nonlocal operator is also nonlinear. In \cite{de2012general} the focus is on a fractional porous medium equation and shows an explicit dependence of the solution on a power of solution, fractional derivative, and initial conditions.

For operators with the structure \eqref{lap}, 
the authors of \cite{burkovska2020affine} show sensitivity with respect to some kernel parameters (such as the size of the support of $\mu$ and its order of singularity $s$). A second-order evolution model inspired from the theory of peridynamics is considered in \cite{coclite2018wellposedness}, where $\mu$ is a nonintegrable kernel ($\mu(x,y)=|y-x|^{-n-2s}$ with $s>0$), and for which the authors show continuity of solutions with respect to initial data. Both papers \cite{burkovska2020affine} and \cite{coclite2018wellposedness} treat the continuous dependence only in the linear setting. To our knowledge, this is the first paper that includes a comprehensive analysis of dependence on boundary data, external forcing, kernel, as well as some certain types of nonlinearities, where the nonlocal interactions are modeled through  weakly singular (i.e. integrable) kernels.   

The stability results of the paper are provided as estimates of the form
\[
    \|u_2-u_1\|_X \leq C\|b_2-b_1\|_Y
\]
for some $C>0$, with appropriate norms $X$ and $Y$, and where $b_i$ denotes data such as the forcing term, boundary data, or the kernel of the nonlocal operator. It is worth noting that simple examples show that, in general, continuity with respect to data may fail even in the setting of local linear elliptic operators. Indeed, consider Hadamard's example \cite{hadamard1952lectures} for the classical Laplace equation. The system:
\[
\begin{cases}
         \displaystyle{u_{xx}+u_{yy}=0}, \quad x\in \re, \, y>0\\
            u(x,0)=0,\ u_y(x,0)=0, \quad x\in \re 
\end{cases}
\]
admits only the solution $v(x,y)=0$. Under a small perturbation in the boundary conditions, however, the system 
\[
\begin{cases}
         \displaystyle{u^n_{xx}+u^n_{yy}=0}, \quad x\in \re,\, y> 0\\
            u^n(x,0)=0, u^n_y(x,0)=\ds\frac{1}{n}\sin(nx), \quad x\in \re, 
\end{cases}
\]
where  $n\in \na$, admits  solutions $u^n(x,y) = \ds\frac{1}{n} e^{-\sqrt{n}} \sin (nx) \sinh (ny)$ of unbounded magnitude as $n\to\infty$. 

In many physical models the nonlocality is exhibited through heterogeneous kernels $\mu(x,y)$. Most commonly, $\mu(x,y)=\tilde{\mu}(|x-y|)$, with $\tilde{\mu}$ a decreasing function with respect to the distance $|x-y|$ between particles. However, if material properties change with the position, one may require an interaction function of the form $\mu(x,y)=\tilde{\mu}(x, |x-y|)$, or of even a more general structure. For example, in geophysics models (see \cite{samko2013fractional} and references within) or medical imaging \cite{yu2015variable}, nonlocal variable fractional operators contain a (nonintegrable) kernel of the form 
\[
\mu(x,y)=\frac{1}{|y-x|^{s(x)}}.
\]
The nonlocality manifests through the kernel, and also through the boundary conditions which, for the well-posedness of the system, must be imposed on sets of positive measure (often referred to as ``collars", when they surround the domain $\Omega$). As experimental data is usually measured only on surfaces, the nonlocal problems raise an additional difficulty as data has to be (artificially) produced on the volumetric boundary; this is related to the ``skin effect", see \cite{ha2011characteristicsdynamicbrittle} for a discussion. Continuous dependence results that quantify the role that variations in boundary conditions have on solutions would alleviate this problem in the nonlocal framework.

The nonlinear setting brings in additional complexities, especially in the nonlocal realm, where classical results (e.g. compactness or embedding theorems, chain rule) are not available.  We are able, however,  to show stability results for $\mc{L}u=f(u)$ for specific nonlinearities $f$, as well as for equations of the type $\mc{L}(h(u))=f$, where $h$ satisfies a lower bound and possesses an invertability property.  

The arguments of the paper rely on: (i) the convolution structure of the integral operators and ensuing properties (such as weighted-mean value formulas, convolution inequalities); (ii) energy-type arguments enabled by the availability of nonlocal versions of integration by parts theorems, as well as Poincar\'e inequalities, and upper bounds for the Poincar\'e constant \cite{foss2019nonlocal}; (iii) estimates that involve the (shrinking) size of the collar $\Gamma$.

 At the numerical level we conduct investigations that validate the bounds obtained theoretically. Different profiles for forcing terms are considered (sinusoidal, sigmoid), boundary data that is discontinuous on the collar, kernels that affect their solutions through different singularities and types of heterogeneities. In the nonlinear case, we perform simulations for varying parameters that control the nonlinearity versus the linear part of the forcing.

\subsection{Contributions of this paper} 

As mentioned above, the paper aims to provide groundwork studies, both theoretical and numerical, in stability of solutions to nonlocal systems. More precisely,
\begin{itemize}
    \item[$\bullet$] We identify exact dependence of solutions on external forcing and boundary data (through  $L^p$ estimates) with two different type of arguments: mean value type theorems and energy estimates.
    \item[$\bullet$] We produce sensitivity studies for nonlocal models with heterogeneous kernels. In comparison to the integration by parts formula as in \cite{andreu2010nonlocal,hinds2012dirichlet}, the heterogeneities bring forth additional terms in integration by parts arguments, for which additional estimates have to be obtained. In particular, one can extract dependence on the horizon size $\delta$ and the degree of (weak) singularity (which were first obtained in \cite{coclite2018wellposedness});
    \item[$\bullet$] For specific nonlinearities we are able to quantify the sensitivity of the solutions to the nonlocal system with respect to the size of the nonlinearity.
    \item[$\bullet$] The numerical studies performed include simulations with discontinuous forcing, discontinuous data on the collar, different types of kernels (symmetric with varying singularity and various heterogeneous ones), as well as nonlinear forcing terms. We investigate the stability of the bounds for the continuous dependence, and the relationship with the theoretical bounds which involve nonlocal Poincar\'e constant (estimated using the arguments of \cite{foss2019nonlocal}. 
\end{itemize}

\subsection{Outline of the paper}  Section 2 of the paper contains preliminaries needed for the proofs, including a list of the main assumptions for the kernel and tools for analysis (such as inequalities and integration by parts). With the background material available, in Section 3 we prove several results on continuous dependence and stability in the linear setting. We consider the nonlinear setting in Section 4 and give various proofs of continuous dependence and stability in the nonlinear setting where we consider nonlinear Laplacian operators, as well as semilinear problems with Lipschitz forcing terms. Section 5 presents numerical studies that validate the results from Sections 3 and 4.

\section{Preliminaries and setup}

\subsection{The setting of nonlocal operators; assumptions and notation}\label{oper}

The results of this paper are set in the framework of nonlocal operators, of which the ones needed are introduced below. The operators are kernel-dependent, which measures the interaction between particles. As in \cite{du2013nonlocal}, for functions $u:\re^n \rightarrow \mathbb{R}$ and $\alpha,\ \mu, \, v:\re^n\times \re^n \rightarrow \mathbb{R}$, we define the nonlocal gradient with kernel $\alpha$ as the two-point operator
\[ \mathcal{G}_\alpha u(x,y):=[u(y)-u(x)]\alpha(x,y), \ x,y\in \re^n. \]
The nonlocal Laplacian with kernel $\mu$ is given by
\[
    \mathcal{L}_{\mu} u (x):=\Int_{\re^n}(u(y)-u(x))\mu(x,y)dy, \quad x \in \re^n.
\]
For symmetric kernels ($\mu(x,y)=\mu(y,x)$) one can write $\mc{L}_{\mu}u=\mc{D}_{\alpha}(\mc{G}_{\alpha}u)$, where the nonlocal divergence of a two-point function is given by
\[
\mc{D}_{\alpha} v(x,y):= \int_{\re^n} v(x,y)\alpha(x,y)-v(y,x)\alpha(y,x) \, dy,
\]
in which case $\mu(x,y)=\alpha^2(x,y)$.

Most results of the paper apply for a large class of heterogeneous  kernels $\mu$ (including anisotropic). In the most basic case we will impose the assumption: 
\begin{list}{\leftmargin=2in \itemindent=-1in}
 \item\textbf{(M1)} $\mu$ is nonnegative and $\mu \in L^1(\re^n\times \re^n)$. 
\end{list}

The domain $\Omega\subseteq\re^n$ (open, bounded set) is surrounded by a collar set $\Gamma$ that has to be chosen appropriately. More precisely, given a kernel $\mu$, we will impose that the collar satisfies the condition: 
\begin{list}{\leftmargin=2in \itemindent=-1in}
 \item\textbf{(M2)}
The domain $\Omega$ (open, bounded set) is surrounded by a collar set $\Gamma$ such that
\[
\bigcup_{x\in\Omega}\text{ supp } \mu(x, \cdot)\setminus \Omega\subseteq\Gamma.
\]
Note that we allow $x\in\dom$ such that $\mu(x,\cdot)$ may not have bounded support.

\end{list}
The coercivity given by the Poincar\'e inequality (see Lemma \ref{Poinc} below) will also need the following lower bound for the kernel in an anulus around the origin. More precisely, 
for each $0<\vep<\delta$, set \[
A_{\vep,\delta}(x):=\{B_\delta(x)\mid |y-x|>\vep\}.
\]

The following are the primary assumptions that we will use.
\begin{list}{\leftmargin=2in \itemindent=-1in}
 \item\textbf{(M3)}  There exists $1\leq p<\infty$, $\mu_0>0$, and $0<\varepsilon<\delta$ such that 
    \begin{equation}\label{kernelgrowth}
    \mu(x,y)\geq \frac{\mu_0}{|y-x|^p},
    \quad \text{for all }x\in\dom\text{ and }
    y\in \text{supp }\mu(x,\cdot)\cap A_{\varepsilon, \delta}(x). 
    \end{equation}
\end{list}
Our last assumptions require integrability for the $x$- and $y$-slices of the kernel.
\begin{list}{\leftmargin=0em \itemindent=0em}
 \item\textbf{(M4)}
 For non-symmetric kernels assume:
  
    \begin{list}{\leftmargin=0in \itemindent=0in}
\item(i)
    For a.e. $y\in\ren$, suppose that $\mu(\cdot,y)\in L^1(\ren)$.
    
    \item(ii) For a.e. $x\in\ren$, suppose that $\mu(x,\cdot)\in L^1(\ren)$.
    \end{list}

    Under this assumption we may introduce the auxiliary functions
\begin{equation}\label{normnote}
    \gamma_\mu(y):=\|\mu(\cdot,y)\|_{L^1(\ren)}
    \nd
     \lambda_\mu(x):=\|\mu(x,\cdot)\|_{L^1(\re^n)}
\end{equation}
and $M_{\mu,p}:=\|\gamma_\mu\|^{1/p}_{L^\infty(\ren)}
\|\lambda_\mu\|^{1-1/p}_{L^\infty(\ren)}$. Notice that if $\mu$ is symmetric then $\gamma_\mu(y)=\lambda_\mu(x)$ for $x=y$ and $M_{\mu,p}=\|\lambda_\mu\|_{L^\infty(\ren)}$.
\end{list}

\vspace*{.2in}
In applications such as peridynamics \cite{Silling}, a prototypical kernel is given by
\begin{equation}\label{mu}
\mu(x,y):=\begin{cases}
\ds\frac{1}{|x-y|^{\beta}}, \quad |x-y|<\delta\\
0, \quad\quad\quad\quad |x-y|\geq \delta,
\end{cases}
\end{equation}
where the parameter $\delta>0$ is called horizon of interaction. For this kernel, the integrability assumption in (M1) is satisfied if $0\leq \beta<n$.

Define the symmetric and antisymmetric parts of $\mu$ by
\begin{equation}\label{symnote}
    \mu_{\sym}(x,y):=\smfrac{1}{2}[\mu(x,y)+\mu(y,x)]
    \nd
    \mu_{\asym}(x,y):=\smfrac{1}{2}[\mu(x,y)-\mu(y,x)].
\end{equation}
We observe that $\mu_{\sym}\ge0$, by assumption (M1). For brevity, we may use 
\begin{equation}\label{symnote2}
\oL_{\sym}=\oL_{\mu_{\sym}} \nd \oG_{\sym}=\oG_{\sqrt{\mu_{\sym}}}. 
\end{equation}
Clearly, $\mu=\mu_{\sym}+\mu_{\asym}$ and $\oL=\oL_{\sym}+\oL_{\asym}$.

\subsection{Tools for analysis}
Some of the proofs below will employ the ``almost" convolution structure of the operator. The results in our paper do not require $\mu(x,y)=\mu(x-y)$, a feature which is amenable to convolution operators. However, a generalization of Young's inequality for this type of more general kernel is available through the following lemma. The following Young's-type inequality is extracted from \cite{okikiolu1970inequalities}.
\begin{lem}\label{L:GenYoung}
Let $1< p< \infty$ be given, and suppose that $\gamma_\mu,\lambda_\mu\in L^\infty(\ren)$. For each $v\in L^p(\dom\cup\bnd)$, define $Tv:\ren\to\ren$ by
\[
    Tv(x)=\int_{\dom\cup\bnd}v(y)\mu(x,y) dy.
\]
Then $Tv\in L^p(\dom)$ and
\begin{equation}\label{Youngs}
    \|Tv\|_{L^p(\dom)}
    \leq
    M_{\mu,p}\|v\|_{L^p(\dom\cup\bnd)}.
\end{equation}
\end{lem}

The nonlocal Laplacian with rotationally symmetric kernel satisfies a list of elliptic-type properties \cite[Prop 3.1]{foss2019bridging}, some of which will be generalized and employed here. The following equality, obtained for convolution kernels in \cite{foss2019bridging}, is a simple consequence of the definition of the nonlocal Laplacian and it will be used in several proofs below:  
\begin{lem} Let $u: \Omega \cup \Gamma \rightarrow \mathbb{R}$ be measurable. Then if $\mu$ satisfies (M1) and $\mathcal{L}_\mu u=f$ in $\Omega$, we have the following property
\begin{equation}\label{eq: MVP}
    u(x)=\frac{1}{\|\mu\|_{L^1(\re^n\times \re^n)}}\Int_{\re^n} u(y)\mu(x,y)dy-\frac{1}{\|\mu\|_{L^1(\re^n\times\re^n)}}f(x), \quad x\in \Omega.
\end{equation}
Note that if $\|\mu\|_{L^1(\re\times\re)}=1$ (e.g. convolution kernels that are probability distributions) and $f(x)=0$, then $u$ satisfies the weighted mean value property:
\[u(x)=\Int_{\re} u(y)\mu(x,y)dy, \quad x\in \Omega.\]
\end{lem}

As key tools for the proofs of our main results we will employ nonlocal versions for integration by parts and a Poincar\'e-type inequality.

\begin{lem}\label{L:IntParts}
Let $u,v:\dom\cup\bnd\to\re$ be measurable. Then
\[
    \int_{\dom\cup\bnd}\oL_{\mu}u(x)v(x)dx
    =
    -\int_{\dom\cup\bnd}\int_{\dom\cup\bnd}
        \oG_{\sym} u (x,y)\oG_{\sym} v (x,y)dydx
    +\int_{\dom\cup\bnd}\int_{\dom\cup\bnd}u(y)v(x)\mu_{\asym}(x,y)dydx
\]
\end{lem}

\begin{proof}
We begin by trivially extending $u,v$ by zero to $\re^n$. Then
\begin{align*}
    \int_{\dom\cup\bnd}\oL_\mu u(x)v(x)dx
    =&
    2\int_{\re^n}\int_{\re^n}[u(y)-u(x)]v(x)\mu_{\sym}(x,y)dydx
    +2\int_{\re^n}\int_{\re^n}[u(y)-u(x)]v(x)\mu_{\asym}(x,y)dydx\\
    =&
    -\int_{\re^n}\int_{\re^n}[u(y)-u(x)][v(y)-v(x)]\mu_{\sym}(x,y)dydx
    +\int_{\re^n}\int_{\re^n}u(y)v(x)\mu_{\asym}dydx.
\end{align*}
Recalling that $u=v=0$ on $\re^n\setminus(\dom\cup\bnd)$ establishes the lemma.
\end{proof}

\begin{lem}\label{L:BndforAsym}
Let  a measurable function $\nu:\re^n\times\re^n\cup\bnd\to\re$. Each of the following holds:
\begin{enumerate}[(a)]
    \item Suppose $\nu$ satisfies both parts of assumption (M4) and that $\gamma_\nu,\lambda_\nu\in L^\infty(\ren)$. Given H\"{o}lder conjugate exponents $1<p,q<\infty$, $u\in L^p(\dom)$, and $v\in L^q(\dom\cup\bnd)$, we have
\[
    \int_{\dom}\int_{\dom\cup\bnd}|u(y)||v(x)||\nu(x,y)|dydx
    \le
    M_{\nu,p}\|u\|_{L^p(\dom)}\|v\|_{L^q(\dom\cup\bnd)}.
\]

    \item Suppose that $\nu\in L^2(\re^n\times\re^n)$. Given $u\in L^2(\dom)$ and $v\in L^2(\dom\cup\bnd)$, we have
\[
    \int_\dom\int_{\dom\cup\bnd}|u(y)||v(x)||\nu(x,y)|dydx
    \le
   \|\nu\|_{L^2(\re^n\times\re^n)}\|u\|_{L^2(\dom)}\|v\|_{L^2(\dom\cup\bnd)}.
\]
\end{enumerate}
\end{lem}
\begin{proof}
Part (a) is a direct consequence of H\"{o}lder's inequality and Lemma~\ref{L:GenYoung}. We write
\begin{align*}
    \int_{\dom}\int_{\dom\cup\bnd}
        |u(x)||v(y)||\nu(x,y)|dydx
    \le&
    \int_{\dom}|u(x)|
        \left(\int_{\dom\cup\bnd}|v(y)|
            |\nu(x,y)dy\right)dx\\
    \le&
    \|u\|_{L^p(\dom)}
    \left(\int_{\dom}\left(
    \int_{\dom\cup\bnd}|v(y)||\nu(x,y)|dy\right)^q dx
        \right)^\frac{1}{q}\\
    \le&
    M_{\nu,q}\|u\|_{L^p(\dom)}\|v\|_{L^q(\dom\cup\bnd)}.
\end{align*}

For part (b), we use Minkowski's integral inequality instead of Lemma~\ref{L:GenYoung}, followed by H\"older's inequality to get
\begin{align*}
    \int_{\dom}\int_{\dom\cup\bnd}|u(x)||v(y)||\nu(x,y)|dydx
    \le&
    \|u\|_{L^2(\dom)}\left(\int_{\dom}\left(
    \int_{\dom\cup\bnd}|v(y)||\nu(x,y)|dy\right)^2dx
        \right)^\frac{1}{2}\\
    &\le
    \|u\|_{L^2(\dom)}\int_{\dom\cup\bnd}\left(
    \int_{\dom}|v(y)|^2
        |\nu(x,y)|^2dx\right)^\frac{1}{2}dy\\
    &\le
    \|\nu(x,y)\|_{L^2(\re^n\times\re^n)}\|u\|_{L^2(\dom)}\|v\|_{L^2(\dom\cup\bnd)}, 
\end{align*}
which gives the conclusion.
\end{proof}

A critical tool for obtaining estimates for solutions to nonlocal problems is the nonlocal Poincar\'e inequality, which can be found in several papers (see for example \cite{aksoylu2011variational}, \cite{andreu2010nonlocal}, \cite{mengesha2014bond}, \cite{ponce2004estimate}). For our results we will need upper bounds for the Poincar\'e constant $C_P$, which can be obtained from \cite[Example 3.2]{foss2019nonlocal}. 

\begin{lem}\label{Poinc}
\cite{foss2019nonlocal} Let $1 \leq p < \infty$, $0  <\varepsilon <\delta$, and an open set $Z \subseteq A(\vep,\delta):=\{x\in \re^n|\, \varepsilon < |x|<\delta \}$ be given. Set $\Gamma:= \cup_{x\in \Omega} \overline{(x+Z)}\setminus\Omega.$ If $u$ is a measurable function over $\Omega$ and $u = 0$ a.e. on $\Gamma$, then
\begin{equation}\Int_\Omega|u(x)|^p dx \leq \frac{\text{diam}(\Omega)^p}{m(Z)}\Int_\Omega \int_Z \frac{|u(x+z)-u(x)|^p}{\|z\|_{\mathbb{R}^n}^p}dz dx,
\end{equation}
where $m(Z)$ is the measure of Z. 

In particular, under assumption (M3) and with $Z:=\text{supp} \,\mu \cap A(\varepsilon, \delta)\,$ we obtain
\begin{equation}\label{Poincineq}
    \Int_\Omega|u(x)|^p dx
    \leq
    C_P\Int_\Omega \Int_{{\text supp}\, \mu} |u(x+z)-u(x)|^p \mu(x,x+z) dz dx,
\end{equation}
where $C_P:=\displaystyle\frac{\text{diam}(\Omega)^p}{\mu_0 m(Z)}$.
\end{lem}

\begin{rem} \label{rem:poincareConstant}
For the numerical results in \S\ref{sec:Numerical}, it will be useful to know the optimal (smallest) Poincar\'e constant $C_P$ in \eqref{Poincineq} for $p=2$ and $\mu(x,y) = 3\delta^{-3}$ on $\Omega = (0,1)$. As $\text{diam}(\Omega)^2=1$, to minimize $C_P$ we must maximize $\mu_0 m(Z)$. From Assumption (M3) we know there exists a $\mu_0>0$ such that 
   \[\frac{3}{\delta^3}=\mu(z)\geq \frac{\mu_0}{|z|^p}, \quad z\in \text{supp }\mu\cap A(\varepsilon, \delta) = A(\varepsilon, \delta). \]
The largest $\mu_0$ results when $z$ is the smallest value in its range, thus $\mu_0 = \varepsilon^2 \mu(z)$ for $p=2$. Note that $m(Z)$ is maximized when the measure of $Z$ is the largest, so $m(Z)=2(\delta-\varepsilon)$. The product $\mu_0 m(Z)$ becomes $3 \delta^{-3}\varepsilon^2 2(\delta-\varepsilon)$, which achieves a maximal value when $\varepsilon = 2 \delta / 3$, and thus the optimal $C_P=9/8$.
\end{rem}

The above conditions guarantee well-posedness of solutions for the linear problems, as well as  some nonlinear problems, as shown in \cite{fossraduwright}.

\section{Continuous dependence of the nonlocal boundary value problem in the linear setting} 
In this section we investigate stability of solutions for the nonlocal Poisson problem
\begin{equation}
    \left\{\begin{array}{ll}
    \displaystyle{\mathcal{L}_{\mu}u(x)=f(x)},
        & \quad x \in \Omega,\\
    u(x)=g(x),
        & \quad x \in \Gamma,
    \end{array}\right.
\end{equation}
under perturbations of the data $f,g$, as well as of the kernel $\mu$. Although the setting is linear for now, some of the methods will extend to the nonlinear setting, which is considered in Section \ref{nonlinearsec}.

\subsection{Stability with respect to the forcing term}

We begin by proving a stability result for solutions under  perturbations of the forcing term by using the mean-value type property.

\begin{theorem}\label{MVPf}
Consider the nonlocal Poisson equations:
\begin{equation}\label{eq:f1}
    \left\{\begin{array}{ll}
    \displaystyle{\mathcal{L}_{\mu}u_i(x)=f_i(x)},
        & \quad x \in \Omega,\\
    u_i(x)=g(x),
        & \quad x \in \Gamma.
    \end{array}\right.
\end{equation}
for $i=1,2$. Given $p \in [1, \infty]$ and $q=\frac{2p}{2p-1}$ and assume the kernel $\mu$ is symmetric and satisfies (M1) and (M2), and in addition, $\mu\in L^q(\mathbb{R}^n \times \mathbb{R}^n)$. Further let $\Omega$ satisfy
\begin{equation}\label{muomega}
m(\Omega)^{\frac{1}{2p}}\frac{\|\mu\|_{L^q(\re^n \times \re^n)}}{\|\mu\|_{L^1(\re^n \times \re^n)}}<1.
\end{equation}
Then, 
\[
\|u_2-u_1\|_{L^{2p}(\Omega)} \leq C_1\|f_2-f_1\|_{L^{2p}(\Omega)}, \]
where the constant $C_1$ above is given by
\begin{equation}\label{const1}
C_1:=\frac{1}
{2\left(\|\mu\|_{L^1(\re^n \times \re^n)}-m(\Omega)^{\frac{1}{2p}}\|\mu\|_{L^q(\re^n \times \re^n)}\right)}.
\end{equation}
\end{theorem}

\begin{proof}
Using Lemma 2.1 and assumptions (M1)--(M2) for $\mu$, we can rewrite \eqref{eq:f1} as
\[u_i(x)=\frac{1}{\|\mu\|_{L^1(\mathbb{R})}}\Int_{\Omega \cup  \Gamma}u_i(y)\mu(x,y)dy-\frac{1}{\|\mu\|_{L^1(\mathbb{R})}}f_i(x).\]
Then, since $\mu(x,y)=\mu(y-x)$, and $u_1=u_2$ in $\re^n\setminus \Omega$ (where we used the fact that $u_1=u_2=g$ on $\Gamma$ and we extended all functions trivially by zero outside $\re^n\setminus \Omega$), we have for all $x\in \re^n$ that

\[|(u_2-u_1)(x)|\leq \frac{1}{\|\mu\|_{L^1(\re^n \times \re^n)}}|(u_2-u_1)*\mu(x)|+\frac{1}{\|\mu\|_{L^1(\re^n \times \re^n)}}|(f_2-f_1)(x)|.\]
Taking the $L^r(\re^n)$ norm of each side and using Young's convolution inequality (with  $\ds\frac{1}{r}+1=\frac{1}{p}+\frac{1}{q}, \, 1\leq p, q, r \leq  \infty$) on the first term, and the fact that $u_2-u_1=f_2-f_1=0$ on $\re^n\setminus \Omega$ we obtain
\begin{equation}
\|u_2-u_1\|_{L^r(\Omega)}\leq 
 \frac{\|\mu\|_{L^q(\re^n \times \re^n)}}{\|\mu\|_{L^1(\re^n \times \re^n)}}\|u_2-u_1\|_{L^p(\Omega)}+\frac{1}{\|\mu\|_{L^1(\re^n \times \re^n)}}\|f_2-f_1\|_{L^r(\Omega)}.
\end{equation}
Now let $r=2p$, so $q=\frac{2p}{2p-1}$. Hence
\begin{equation}\label{ineq1}
\|u_2-u_1\|_{L^{2p}(\Omega)}\leq \frac{\|\mu\|_{L^q(\re^n \times \re^n)}}{\|\mu\|_{L^1(\re^n \times \re^n)}} \|u_2-u_1\|_{L^p(\Omega)}+\frac{1}{\|\mu\|_{L^1(\re^n \times \re^n)}} \|f_2-f_1\|_{L^{2p}(\Omega)}.
\end{equation}
From H\"older's inequality we have that
\begin{equation*}
\|u_2-u_1\|_{L^p(\Omega)} \leq m(\Omega)^{\frac{1}{2p}}\|u_2-u_1\|_{L^{2p}(\Omega)} ,
\end{equation*}
so from \eqref{ineq1} we obtain
\[\|u_2-u_1\|_{L^{2p}(\Omega)}\leq m(\Omega)^{\frac{1}{2p}}\frac{\|\mu\|_{L^q(\re^n \times \re^n)}}{\|\mu\|_{L^1(\re^n \times \re^n)}} \|u_2-u_1\|_{L^{2p}(\Omega)}+\frac{1}{\|\mu\|_{L^1(\re^n \times \re^n)}} \|f_2-f_1\|_{L^{2p}(\Omega)}.\]
Thus, under the assumption $m(\Omega)^{\frac{1}{2p}}\frac{\|\mu\|_{L^q(\re^n \times \re^n)}}{\|\mu\|_{L^1(\re^n \times \re^n)}}<1$ we obtain
\[
\|u_2-u_1\|_{L^{2p}(\Omega)}\leq
\frac{1}
{\|\mu\|_{L^1(\re^n \times \re^n)}-m(\Omega)^{\frac{1}{2p}}\|\mu\|_{L^q(\re^n \times \re^n)}}
\|f_2-f_1\|_{L^{2p}(\Omega)}.
\]
\end{proof}

\begin{rem}
Notice that this result usually requires a large collar in order for the condition  \eqref{muomega} to hold. Indeed, for a typical peridynamic kernel as given by \eqref{mu}, the condition \eqref{muomega} becomes
\[
m(\Omega)<\delta^{\beta\frac{2p}{2p-1}},
\]
so for constant kernels ($\beta=0$) we need to impose $m(\Omega)<1.$ However, Theorem \ref{MVPf} does yield stability results for all $L^p$ norms with $p\geq 2$. Additionally, as we will see in the sequel, the proof generalizes to certain nonlinear problems (see Theorem \ref{nonlinear1}).  Next, with a similar argument, we establish an alternative to the stability result above that replaces the $m(\Omega)$ constraint with a restriction on $\mu$ and allows for an asymmetric component to the kernel. Alternatively, we can obtain a similar stability result (but only in $L^2$) by using an energy argument with no requirement on the size of the domain as in Theorem \ref{linearforcingenergy}.
\end{rem}

\begin{theorem}\label{MVPfasym}
Consider the nonlocal Poisson equations:
\begin{equation}\label{eq:f1asym}
    \left\{\begin{array}{ll}
    \displaystyle{\mathcal{L}_{\mu}u_i(x)=f_i(x)},
        & \quad x \in \Omega,\\
    u_i(x)=g(x),
        & \quad x \in \Gamma.
    \end{array}\right.
\end{equation}
for $i=1,2$. Let $r \geq 1$. If (M1), (M2), and (M4) hold, we have
\[\|u_2-u_1\|_{L^{r}(\Omega)} \leq C_2\|f_2-f_1\|_{L^{r}(\Omega)}. \]
The constant $C_2$ above is given by
\begin{equation}\label{const2}
C_2:=\frac{\|\frac{1}{\lambda}\|_{L^\infty(\Omega \cup \Gamma)}}{1- M_{\mu,r}\|\frac{1}{\lambda}\|_{L^{\infty}(\Omega \cup \Gamma)}}.
\end{equation}
\end{theorem}

\begin{proof}
 From \eqref{eq:f1asym} we have 
\[u_i(x)=\frac{1}{\lambda(x)}\int_{\Omega \cup \Gamma} u_i(y)\mu(x,y)dy-\frac{1}{\lambda(x)}f_i(x)=\frac{1}{\lambda(x)} \int_{\Omega\cup \Gamma} u_i(y)\mu(x,y)dy-\frac{1}{\lambda(x)}f_i(x)\].

\noindent Subtracting the two solutions and taking the $L^r$ norm, we have 
\begin{equation}
\begin{split}
 &\|u_2-u_1\|_{L^r(\Omega \cup \Gamma)}\\
  \leq &  \underbrace{\left(\int_{\Omega\cup \Gamma}\left|\frac{1}{\lambda(x)}\left(\int_{\Omega\cup \Gamma}\left| (u_2-u_1)(y)\mu(x,y)\right|dy\right)\right|^r dx\right)^{1/r}}_{=:I}+\underbrace{\left(\int_{\Omega\cup \Gamma}\left|\frac{1}{\lambda(x)}(f_2-f_1)\right|^r dx\right)^{1/r}}_{=:II}.\\
  \end{split}
\end{equation}

\noindent We handle $I$ first. By H\"older's inequality with $p,q\geq1$ and $\frac{1}{p}+\frac{1}{q}=1$
we have 
\begin{equation}\label{estI}
\begin{split}
  &\left(\int_{\Omega\cup \Gamma}\left|\frac{1}{\lambda(x)}\left(\int_{\Omega\cup \Gamma}\left| (u_2-u_1)(y)\mu(x,y)\right|dy\right)\right|^r dx\right)^{1/r}\\
  \leq& \left(\int_{\Omega \cup \Gamma}\frac{1}{\left|\lambda(x)\right|^{pr}}dx\right)^{1/pr}\left( \int_{\Omega \cup \Gamma} \left|\int_{\Omega\cup\Gamma}\left| (u_2-u_1)(y)\mu(x,y)\right|dy\right|^{rq}dx\right)^{1/rq}.
  \end{split}
\end{equation}
By Lemma \ref{L:GenYoung}, we have
\begin{align*}
  & \left(\int_{\Omega \cup \Gamma}\left| \int_{\Omega\cup\Gamma}\left| (u_2-u_1)(y)\mu(x,y)\right|dy\right|^{rq}dx\right)^{1/rq}\leq  M_{\mu,rq} \|u_2-u_1\|_{L^{rq}(\Omega \cup \Gamma)}.\\
\end{align*}
Thus from \eqref{estI}, by letting $q=1$ (so that $rq=r$ and $p=\infty$) and denoting by $s$ the H\"older conjugate of $r$ ($\frac{1}{r}+\frac{1}{s}=1$) we have 
\begin{align*}
   \left(\int_{\Omega\cup \Gamma}\left|\frac{1}{\lambda(x)}\left(\int_{\Omega\cup \Gamma}\left| (u_2-u_1)(y)\mu(x,y)\right|dy\right)\right|^r dx\right)^{1/r} \leq  
  M_{\mu,r}\left\|\frac{1}{\lambda}\right\|_{L^{\infty}(\Omega \cup \Gamma)}\|u_2-u_1\|_{L^{r}(\Omega\cup \Gamma)}.
\end{align*}

\noindent Next we handle $II$. By H\"older's 
inequality and letting $p=\infty$ and $q=1$, we have 
\begin{align*}
    \int_{\Omega \cup \Gamma}\left|\frac{1}{\lambda(x)}(f_2-f_1)(x)\right|dx\leq  
\left\|\frac{1}{\lambda}\right\|_{L^\infty(\Omega \cup \Gamma)} \|f_2-f_1\|_{L^{1}(\Omega \cup \Gamma)}.
\end{align*}
\noindent Thus we have 
\[\|u_2-u_1\|_{L^r(\Omega \cup \Gamma)}\leq  M_{\mu,r}\left\|\frac{1}{\lambda}\right\|_{L^{\infty}(\Omega \cup \Gamma)}\|u_2-u_1\|_{L^{r}(\Omega\cup \Gamma)}+ \left\|\frac{1}{\lambda}\right\|_{L^\infty(\Omega \cup \Gamma)} \|f_2-f_1\|_{L^{1}(\Omega \cup \Gamma)},
\]
which we can rewrite as 
\[\|u_2-u_1\|_{L^r(\Omega \cup \Gamma)}\leq \frac{\|\frac{1}{\lambda}\|_{L^\infty(\Omega \cup \Gamma)}}{1- M_{\mu,r}\|\frac{1}{\lambda}\|_{L^{\infty}(\Omega \cup \Gamma)}}\|f_2-f_1\|_{L^1(\Omega \cup \Gamma)}.\]
\end{proof}

\begin{rem} Notice that the restriction $1- M_{\mu,r}\|\frac{1}{\lambda}\|_{L^{\infty}(\Omega \cup \Gamma)}>0$ can not be satisfied by kernels $\mu$ radially symmetric and positive, as:
\begin{align*}
\|\gamma\|^{1/p}_{L^\infty(\ren)}
\|\lambda\|^{1-1/p}_{L^\infty(\ren)}\left\|\frac{1}{\lambda}\right\|_{L^{\infty}(\Omega \cup \Gamma)}=\|\lambda\|_{L^\infty(\ren)}\left\|\frac{1}{\lambda}\right\|_{L^{\infty}(\Omega \cup \Gamma)}\geq 1.
\end{align*}

However, one can generate simple examples of kernels that  fit this restriction. For example, let $p=2$, $n=1$, \, $0<a<b$ and assume that $\Omega=(a,b).$  Let $\mu(x,y)=x$ on $B_\delta(x)$. Then we have 
\begin{align*}
    &\|\gamma\|^{1/2}_{L^\infty(\ren)}
\|\lambda\|^{1/2}_{L^\infty(\ren)}\left\|\frac{1}{\lambda}\right\|_{L^{\infty}(\Omega \cup \Gamma)}\\
&=\left(\max_{y\in (x-\delta, x+\delta)} \int_a^b|x|dx\right)^{\frac{1}{2}}  \left( \max_{x\in(a,b)} \int_{x-\delta}^{x+\delta}|x|dy\right)^{\frac{1}{2}}  \left( \max_{x\in(a,b)} \frac{1}{\int_{x-\delta}^{x+\delta}|x|dy} \right)\\
    =&\left(\max_{y\in(x-\delta, x+\delta)} \left(\frac{b^2}{2}-\frac{a^2}{2}\right)\right)^{\frac{1}{2}}  \left( \max_{x\in(a,b)} 2\delta x\right)^{\frac{1}{2}}  \left( \max_{x\in(a,b)} \frac{1}{2\delta x} \right)\\
    =&\left(\frac{b^2}{2}-\frac{a^2}{2}\right)^{\frac{1}{2}}  \left( 2\delta b\right)^{\frac{1}{2}}  \left( \frac{1}{2\delta a} \right)=\frac{(b^2-a^2)^{\frac{1}{2}}}{2} \left( \frac{ b^{\frac{1}{2}} }{\delta^{\frac{1}{2}} a} \right)<1.
\end{align*}
There are a variety of choices for $a,b,\delta$ that satisfy this restriction. Indeed, take for example, $a=1$, then for $b<(-1+\sqrt{17})/2$, we have that there exist $\delta>0$ that satisfy the restriction above as well as $\delta<m(\Omega)=b-a$.

Notice that if instead, we let $\mu(x,y)=y$, we have a slightly different restriction since:
\begin{align*}
 &\|\gamma\|^{1/2}_{L^\infty(\ren)}
\|\lambda\|^{1/2}_{L^\infty(\ren)}\left\|\frac{1}{\lambda}\right\|_{L^{\infty}(\Omega \cup \Gamma)}\\
    &=\left(\max_{y\in(x-\delta, x+\delta)} \int_a^b|y|dx\right)^{\frac{1}{2}}  \left(\max_{x\in(a,b)} \int_{x-\delta}^{x+\delta}|y|dy\right)^{\frac{1}{2}}  \left( \max_{x\in(a,b)} \frac{1}{\int_{x-\delta}^{x+\delta}|y|dy} \right)\\
    =&\left(\max_{y\in(x-\delta, x+\delta)} \left(y(b-a)\right)\right)^{\frac{1}{2}}  \left( \max_{x\in(a,b)} \delta x\right)^{\frac{1}{2}}  \left( \max_{x\in(a,b)} \frac{1}{\delta x} \right)\\
    =&\left( (x+\delta)(b-a) \right)^{\frac{1}{2}}  \left(\delta b\right)^{\frac{1}{2}}  \left( \frac{1}{\delta a} \right)=\left( (x+\delta)(b-a) \right)^{\frac{1}{2}} \left( \frac{b^{\frac{1}{2}} }{\delta^{\frac{1}{2}} a} \right)<1,
\end{align*}
where $a\neq 0$. The restriction above must hold for all $x \in (a,b)$, so we impose 
\begin{align*}
   \left((b+\delta)(b-a) \right)^{\frac{1}{2}} \left( \frac{b^{\frac{1}{2}} }{\delta^{\frac{1}{2}} a} \right)<1,
\end{align*}
or equivalently,
\begin{equation}\label{restrict1}
   b^3-ab^2<\delta(ab+a^2-b^2).
\end{equation}
Notice that since the left hand side is always positive, we need $ab+a^2-b^2>0$, which implies $a<b<\frac{1+\sqrt{5}}{2}a$, thus severely restricting the length of the interval $(a,b).$ Additionally, by letting $b=\alpha a$ for $\alpha>1$, then \eqref{restrict1} becomes $\delta>\frac{a\alpha^3(\alpha-1)}{\alpha+1-\alpha^2}$. A simple calculation shows that this inequality implies $\delta>a(\alpha-1)$, which is the length of $(a,b)$. Thus, no $\delta$ exists such that $\delta<b-a$, so the collar size exceeds the size of the domain, a similar restriction to that in Theorem \ref{MVPf}.
\end{rem}

\begin{theorem}\label{linearforcingenergy}
Consider the nonlocal Poisson's equation over the domain $\Omega\subset \mathbb{R}$.
\begin{equation}
    \left\{\begin{array}{ll}
    \displaystyle{\mathcal{L}_{\mu}u_i(x)=f_i(x)},
        & \quad x \in \Omega,\\
    u_i(x)=g(x),
        & \quad x \in \Gamma.
    \end{array}\right.
\end{equation}
for $i=1,2$. Let $g\in L^2(\Gamma)$.
Then, if (M3) and (M4) are satisfied and $1>M_{\mu_{\text{asym}},2}C_P$,
\[ \| u_2-u_1\|_{L^2(\Omega)}
    \leq\frac{C_P}{1-M_{\mu_{\text{asym}},2}C_P}\|f_2-f_1\|_{L^2(\Omega)}\]
where $C_P$ is the Poincaré constant from \cite{foss2019nonlocal}.
\end{theorem}

\begin{proof}
 Multiplying $\mathcal{L}_{\mu_2}(u_2-u_1)$ by $u_2-u_1$, integrating, and using the notation of \eqref{symnote2} we have
\begin{align*}
    \int_{\Omega}(u_2(x)-u_1(x))\mathcal{L}_{\mu}(u_2-u_1)(x) dx &=-\int_{\Omega}\int_{\Omega}
    [\mathcal{G}_{\text{sym}}(u_2-u_1)]^2 dy dx\\
    &\hspace{20pt}
    +\int_{\Omega}\int_{\Omega}[u_2(y)-u_1(y)][u_2(x)-u_1(x)]
    \mu_{\text{asym}}(x,y)dydx.
\end{align*}
By H\"older's inequality we have
\begin{equation} \label{asymb}
\begin{split}
&\left|\int_{\Omega}\int_{\Omega}
        [u_2(y)-u_1(y)][u_2(x)-u_1(x)]\mu_{\text{asym}}(x,y)dydx\right|\\ 
    &\leq\|u_2-u_1\|_{L^2(\Omega)}
    \left(\int_{\Omega}\left(
    \int_{\Omega}|u_2(y)-u_1(y)||\mu_{\text{asym}}(x,y)|dy\right)^2dx
        \right)^\frac{1}{2}\\
    &\leq M_{\mu_{\text{asym}},2}\|u_2-u_1\|_{L^2(\Omega)}^2.
\end{split}
\end{equation}
For the last line, Lemma \ref{L:GenYoung} was used and $\gamma_{2,\text{asym}}\in L^\infty(\Omega\cup\Gamma)$ is defined by $\gamma_{\text{asym}}(y):=\|\mu_{\text{asym}}(\cdot,y)\|_{L^1(\mathbb{R}^n)}$.

Then from nonlocal Poincar\'e's inequality, we have that
\begin{equation} \label{asymb2}
\begin{split}
   \| u_2-u_1\|_{L^2(\Omega)}^2
   &\leq C_P \|\mathcal{G}_{\text{sym}}(u_2-u_1)\|_{L^2(\Omega)}^2 
   =C_P\int_{\Omega}\int_{\Omega}
   [\mathcal{G}_{\text{sym}}(u_2-u_1)]^2 dy dx \\
   &\leq\left|
   C_P\int_{\Omega}\int_{\Omega}[u_2(y)-u_1(y)][u_2(x)-u_1(x)]
    \mu_{\text{asym}}(x,y)dydx\right|\\
    &\hspace{50pt}+\left|
    C_P\int_{\Omega \cup \Gamma}(u_2(x)-u_1(x))\mathcal{L}_{\mu}(u_2-u_1)(x) dx\right|\\
    &\le
    C_PM_{\mu_{\text{asym}},2}\|u_2-u_1\|_{L^2(\Omega)}^2
    +C_P\left|\int_{\Omega}(u_2(x)-u_1(x))
        \mathcal{L}_{\mu}(u_2-u_1)(x)dx\right|.
\end{split}
\end{equation}

The second term is bounded above (using Hölder's inequality) by $C_P\|u_2-u_1\|_{L^2(\Omega)}\|f_2-f_1\|_{L^2(\Omega)}.$

And so 
\begin{align*}
   \| u_2-u_1\|_{L^2(\Omega)}^2
    &\leq
    C_PM_{\mu_{\text{asym}},2}\|u_2-u_1\|_{L^2(\Omega)}^2
    +C_P\|u_2-u_1\|_{L^2(\Omega)}\|f_2-f_1\|_{L^2(\Omega)}.
\end{align*}

Rearranging, we have 
\begin{align*}
   \| u_2-u_1\|_{L^2(\Omega)}
    &\leq\frac{C_P}{1-M_{\mu_{\text{asym}},2}C_P}\|f_2-f_1\|_{L^2(\Omega)}.
\end{align*}
\end{proof}


\subsection{Continuous Dependence on the Collar Term}

For the continuous dependence on the collar, we only consider the energy argument. Following the ``mean value" type argument yields a condition equivalent to needing support for the kernel to be larger than both the domain and the collar, rendering the result useless.

\begin{corollary}
Consider the nonlocal Poisson's equation over the domain $\Omega\subset \mathbb{R}$.
\begin{equation}\label{eqcollar}
    \left\{\begin{array}{ll}
    \displaystyle{\mathcal{L}_{\mu}u_i(x)=f(x)},
        & \quad x \in \Omega,\\
    u_i(x)=g_i(x),
        & \quad x \in \Gamma.
    \end{array}\right.
\end{equation}
for $i=1,2$. Let $g\in L^2(\Gamma)$.
Then, if (M3) and (M4) are satisfied and $1>M_{\mu_{\text{asym}},2}C_P$,
\[\| u_2-u_1\|_{L^2(\Omega \cup \Gamma)} \leq \frac{C_P}{1-M_{\mu_{\text{asym}},2}C_P} \|\mathcal{L}_{\mu}(g_2-g_1)\|_{L^2(\Omega \cup \Gamma)}\]
where $C_P$ is the Poincaré constant from \cite{foss2019nonlocal}.
\end{corollary}

\begin{proof}
We begin by considering $w_i=u_i-g_i$. We extend $g_i$ by $0$ to $\Omega$, since $g\in L^2(\Gamma)$. Then by linearity of the nonlocal Laplacian \eqref{eqcollar} with the function $w_i$ is
\begin{equation*}
    \left\{\begin{array}{ll}
    \displaystyle{\mathcal{L}_{\mu}w_i=f-\mathcal{L}_{\mu}g_i},
        & \quad x \in \Omega,\\
    w_i=0,
        & \quad x \in \Gamma.
    \end{array}\right.
\end{equation*}
From here, we apply Theorem 3.3 to have 
\begin{align*}
    \| u_2-u_1\|_{L^2(\Omega \cup \Gamma)} \leq \frac{C_P}{1-M_{\mu_{\text{asym}},2}C_P}\|\mathcal{L}_{\mu}g_2-\mathcal{L}_{\mu}g_1\|_{L^2(\Omega\cup \Gamma)}.
\end{align*}
Then since the support of $g$ is $\Gamma$, we yield 
\begin{align*}
    \| u_2-u_1\|_{L^2(\Omega \cup \Gamma)} \leq\frac{C_P}{1-M_{\mu_{\text{asym}},2}C_P}\|\mathcal{L}_{\mu}(g_2-g_1)\|_{L^2(\Gamma)}.
\end{align*}
\end{proof}

\begin{theorem}\label{collarenergy}
Consider the nonlocal Poisson's equation over the domain $\Omega\subset \mathbb{R}$.
\begin{equation}
    \left\{\begin{array}{ll}
    \displaystyle{\mathcal{L}_{\mu}u_i(x)=f(x)},
        & \quad x \in \Omega,\\
    u_i(x)=g_i(x),
        & \quad x \in \Gamma.
    \end{array}\right.
\end{equation}
for $i=1,2$. Then, if (M3) and (M4) are satisfied and $1>C_PM_{\mu_{\text{asym}},2}$,
\begin{align*}
     \| u_2-u_1\|_{L^2(\Omega)}\leq 
   C\|g_2-g_1\|_{L^2(\Gamma)},\\
\end{align*}
where $C=\frac{C_P\|\mu\|_{L^2(\Omega\times \Gamma)}}{1-C_PM_{\mu_{\text{asym}},2}}$.
\end{theorem}

\begin{proof}
Define $v:=u_1$ in $\Omega$ and $v:=u_2$ in $\Gamma$.
Notice that 
\begin{align*}
     \int_{\Omega\cup \Gamma}(u_2(x)-v(x))\mathcal{L}_{\mu}(u_2-v)(x) dx&= \int_{\Omega}(u_2(x)-u_1(x))\mathcal{L}_{\mu}(u_2-v)(x) dx\\
     &= \int_{\Omega}(u_2(x)-u_1(x))\int_{\Gamma}(g_2(y)-g_1(y))\mu(x,y)dy dx\\
      &\leq\|g_2-g_1\|_{L^2(\Gamma)}\|\mu\|_{L^2(\Omega\times \Gamma)}\|u_2-u_1\|_{L^2(\Omega)}\\
\end{align*}
Multiplying $\mathcal{L}_{\mu_2}(u_2-v)$ by $u_2-v$ and integrating, we have
\begin{align*}
    \int_{\Omega\cup \Gamma}(u_2(x)-v(x))\mathcal{L}_{\mu}(u_2-v)(x) dx &=-\int_{\Omega\cup \Gamma}\int_{\Omega\cup \Gamma}
    [\mathcal{G}_{\text{sym}}(u_2-v)]^2 dy dx\\
    &\hspace{20pt}
    +\int_{\Omega\cup \Gamma}\int_{\Omega\cup \Gamma}[u_2(y)-v(y)][u_2(x)-v(x)]
    \mu_{\text{asym}}(x,y)dydx.
\end{align*}
Here $\mathcal{G}_{\text{sym}}=\mathcal{G}_{\sqrt{\mu_{\sym}}}$. Similar to \eqref{asymb}, we have
\begin{align*}
    &\int_{\Omega\cup \Gamma}\int_{\Omega\cup \Gamma}[u_2(y)-v(y)][u_2(x)-v(x)]
    \mu_{\text{asym}}(x,y)dydx\\
    &\leq \left|\int_{\Omega}\int_{\Omega}
        [u_2(y)-u_1(y)][u_2(x)-u_1(x)]\mu_{\text{asym}}(x,y)dydx\right|\\
        &\leq M_{\mu_{\text{asym}},2}\|u_2-u_1\|_{L^2(\Omega)}^2.
\end{align*}
Then from nonlocal Poincar\'e's inequality (since $u_2-v=0$ on $\Gamma$), 
we have that
\begin{align*}
   \| u_2-u_1\|_{L^2(\Omega)}^2&= \| u_2-v\|_{L^2(\Omega)}^2\\
   &\leq 
   C_P\int_{\Omega}\int_{\Omega}[u_2(y)-v(y)][u_2(x)-v(x)]
    \mu_{\text{asym}}(x,y)dydx+   C_P\int_{\Omega\cup \Gamma}(u_2(x)-v(x))\mathcal{L}_{\mu}(u_2-v) dx\\
       &\le
    C_PM_{\mu_{\text{asym}},2}\|u_2-u_1\|_{L^2(\Omega)}^2+ C_P\|g_2-g_1\|_{L^2(\Gamma)}\|\mu\|_{L^2(\Omega\times \Gamma)}\|u_2-u_1\|_{L^2(\Omega)}.
\end{align*}
Consequently, 
\begin{align*}
   \| u_2-u_1\|_{L^2(\Omega)}^2\leq 
    C_PM_{\mu_{\text{asym}},2}\|u_2-u_1\|_{L^2(\Omega)}^2+ C_P\|g_2-g_1\|_{L^2(\Gamma)}\|\mu\|_{L^2(\Omega\times \Gamma)}\|u_2-u_1\|_{L^2(\Omega)}.
\end{align*}
Rearranging, we have 
\begin{equation*}
    \| u_2-u_1\|_{L^2(\Omega)}\leq 
   \frac{C_P\|\mu\|_{L^2(\Omega\times \Gamma)}}{ 1-C_PM_{\mu_{\text{asym}},2}}\|g_2-g_1\|_{L^2(\Gamma)}.
\end{equation*}
\end{proof}

\subsection{Stability with Changes in the Kernel}
Lastly we consider the stability of the solution due to perturbations of the kernel, which alters the operator itself. From the following result, we have both the $L^2$ and the $L^\infty$ norms.

To provide a concise statement for our next theorem, we introduce some supplementary notation. Recall the definition for $\lambda_{\mu_i}$ in~\eqref{normnote}. For convenience, we will use $\lambda_i=\lambda_{\mu_i}$. Define the normalized kernels $\wt{\mu}_i\in L^1(\re^n\times\re^n)$ by 
\begin{equation}\label{normkernel}
    \wt{\mu}_i(x,y):=\lc\begin{array}{ll}
        \lambda_i(x)^{-1}\mu_i(x,y), & x\in\dom\\
        0, & x\in\ren\bs\dom.
    \end{array}\rt
\end{equation}
(Note that, by assumption (M3), we find $\lambda_i(x)>0$ for each $x\in\dom$). We will also write
\[
    \wt{\gamma}_i=\gamma_{\wt{\mu}_i},\quad
    \wt{\gamma}_{1,2}=\gamma_{\wt{\mu}_2-\wt{\mu}_1},\quad
    \wt{\gamma}_{i,\asym}=\gamma_{\wt{\mu}_{i,\asym}},\nd
    \wt{\lambda}_{i,\asym}=\lambda_{\wt{\mu}_{i,\asym}}.
\]

\begin{theorem} \label{thm:kernelBound}
Consider the nonlocal Poisson's equation over the domain $\Omega\subset \mathbb{R}$.
\begin{equation}\label{E:Poisson}
    \left\{\begin{array}{ll}
    \displaystyle{\mathcal{L}_{\mu_i}u_i(x)=f(x)},
        & \quad x \in \Omega,\\
    u_i(x)=g(x),
        & \quad x \in \Gamma.
    \end{array}\right.
\end{equation}
for $i=1,2$. Define $K=\lp\frac{1}{\lambda_2}-\frac{1}{\lambda_1}\rp$. Suppose that $\mu_i$ satisfies (M3) and (M4).
\begin{enumerate}[(a)]
\item If $\wt{M}_2:=\|\wt{\gamma}_{2,\asym}\|_{L^\infty(\ren)}\|\wt{\lambda}_{2,\asym}\|_{L^\infty(\ren)}<C_P^{-1}$, then
\[
    \|u_2-u_1\|_{L^2(\dom)}\le\frac{C_P}{1-C_P\wt{M}_2}
    \ls2\|\wt{\gamma}_{1,2}\|_{L^\infty(\ren)}
    \|\wt{\lambda}_{2}-\wt{\lambda}_{1}\|_{L^\infty(\ren)}\|u_1\|_{L^2(\dom\cup\bnd)}
    +
   \|K\|_{L^\infty(\dom)}\|f\|_{L^2(\dom)}\rs.
\]
\item If $\mu_i\in L^2(\ren\times\ren)$ and $\|\wt{\mu}_{2,\asym}\|_{L^2(\ren\times\ren)}<C_P^{-1}$, then
\[
    \|u_2-u_1\|_{L^2(\dom)}\le\frac{C_P}{1-C_P\|\wt{\mu}_{2,\asym}\|_{L^2(\ren\times\ren)}}
    \ls2\|\wt{\mu}_2-\wt{\mu}_1\|_{L^2(\ren\times\ren)}\|u_1\|_{L^2(\dom\cup\bnd)}
    +
   \|K\|_{L^\infty(\dom)}\|f\|_{L^2(\dom)}\rs.
\]
\end{enumerate}
\end{theorem}

\begin{proof}
Note $K\in L^\infty(\dom)$ by (M4) and the normalized kernels $\wt{\mu_i}, \, i=1,2$ satisfy (M4).
First, we establish a simple identity that follows from the Poisson equation. We can rewrite~\eqref{E:Poisson} as
\begin{align*}
    u_i(x)
    &=
    \Int_{\Omega \cup  \Gamma}u_i(y)\widetilde{\mu_i}(x,y)dy-\frac{f(x)}{\lambda_i(x)}.
\end{align*}
Using the definition of $K$, we may write
\begin{align*}
    u_2(x)-u_1(x)&=\Int_{\Omega \cup  \Gamma}u_2(y)\widetilde{\mu}_2(x,y)dy-\Int_{\Omega \cup  \Gamma}u_1(y)\widetilde{\mu}_1(x,y)dy-K(x)f(x)\\
    &=\Int_{\Omega}(u_2(y)-u_1(y)\widetilde{\mu}_2(x,y)dy
    -\Int_{\Omega \cup  \Gamma}u_1(y)(\widetilde{\mu}_1(x,y)-\widetilde{\mu}_2(x,y))dy
    -K(x)f(x).
\end{align*}
Rearranging, we obtain
\begin{align*}
\nonumber
   \Int_{\Omega \cup  \Gamma}u_1(y)(\widetilde{\mu}_1(x,y)-\widetilde{\mu}_2(x,y))dy
   &=\Int_{\Omega}(u_2(y)-u_1(y)\widetilde{\mu}_2(x,y)dy- (u_2(x)-u_1(x)+K(x)f(x))\\
   &=\mathcal{L}_{\widetilde{\mu}_2}(u_2-u_1)-K(x)f(x),
\end{align*}
and thus
\begin{equation}\label{E:OpConv}
    \mathcal{L}_{\widetilde{\mu}_2}(u_2-u_1)(x)
    =
    \int_{\dom\cup\bnd}
    u_1(y)(\widetilde{\mu}_1(x,y)-\widetilde{\mu}_2(x,y))dy
    +K(x)f(x).
\end{equation}

We next employ the nonlocal Poincar\'{e} inequality to bound $\|u_2-u_1\|_{L^2(\dom)}$. Multiplying $\mathcal{L}_{\widetilde{\mu}_2}(u_2-u_1)$ with $u_2-u_1$ and using Lemma~\ref{L:IntParts} produces
\begin{align*}
    \Int_{\dom\cup\bnd}(u_2(x)-u_1(x))\mathcal{L}_{\widetilde{\mu}_2}(u_2-u_1)(x) dx &=-\Int_{\Omega\cup\bnd}\Int_{\Omega\cup\bnd}
    \ls\wt{\oG}_{2,\sym}(u_2-u_1)(x)\rs^2 dy dx\\
    &\hspace{20pt}
    +\int_{\dom\cup\bnd}\int_{\dom\cup\bnd}[u_2(y)-u_1(y)][u_2(x)-u_1(x)]
    \wt{\mu}_{2,\asym}(x,y)dydx.
\end{align*}
Here we used $\wt{\oG}_{2,\sym}=\oG_{\sqrt{\wt{\mu}_{2,\sym}}}$ and the fact that $u_2-u_1=0$ on $\bnd$. Rearranging and using Lemma~\ref{Poinc}, we obtain
\begin{align}
\label{E:MainIneq}
   \| u_2-u_1\|_{L^2(\Omega)}^2
   &\leq C_P \|\wt{\oG}_{2,\sym}(u_2-u_1)\|_{L^2(\Omega)}^2 
   =C_P\Int_{\Omega}\Int_{\Omega}
   \wt{\oG}_{2,\sym}(u_2-u_1)]^2 dy dx \\
\nonumber
   &=
   C_P\underbrace{\int_{\dom}\int_{\dom}[u_2(y)-u_1(y)][u_2(x)-u_1(x)]
    \wt{\mu}_{2,\asym}(x,y)dydx}_{=:I}\\
\nonumber
    &\hspace{50pt}-
    C_P\underbrace{\Int_{\Omega}
    (u_2(x)-u_1(x))\mathcal{L}_{\widetilde{\mu}_{2}}(u_2-u_1)(x) dx}_{=:II}.
\end{align}
For $I$, Lemma~\ref{L:BndforAsym}(a) provides
\begin{equation}\label{E:BndIa}
    |I|\le \|\wt{\gamma}_{2,\asym}\|_{L^\infty(\ren)}
    \|\wt{\lambda}_{2,\asym}\|_{L^\infty(\ren)}\|u_2-u_1\|^2_{L^2(\dom)}
    =\wt{M}_2\|u_2-u_1\|^2_{L^2(\dom)},
\end{equation}
while Lemma~\ref{L:BndforAsym}(b) gives us 
\begin{equation}\label{E:BndIb}
    |I|\le \|\wt{\mu}_{2,\asym}\|_{L^2(\ren\times\ren)}\|u_2-u_1\|^2_{L^2(\dom)}.
\end{equation}
Note that  $\wt{\mu_i}\in L^2(\ren\times\ren)$ by assumption (b) of the Theorem.

For $II$, we use~\eqref{E:OpConv} and get
\begin{align*}
    |II|
    \le&
    \int_{\dom}\int_{\dom\cup\bnd}
    |u_2(x)-u_1(x)||u_1(y)||\widetilde{\mu}_1(x,y)-\widetilde{\mu}_2(x,y)|dydx
    +\int_\dom|u_2(x)-u_1(x)||K(x)||f(x)|dx\\
    \le&
    \int_{\dom}\int_{\dom\cup\bnd}
    |u_2(x)-u_1(x)||u_1(y)||\widetilde{\mu}_1(x,y)-\widetilde{\mu}_2(x,y)|dydx
    +\|K\|_{L^\infty(\dom)}\|f\|_{L^2(\dom)}\|u_2-u_1\|_{L^2(\dom)}.
\end{align*}
Similarly Lemma~\ref{L:BndforAsym} yields either
\begin{equation}\label{E:BndIIa}
    |II|\le \|\wt{\gamma}_{1,2}\|_{L^\infty(\ren)}
    \|\wt{\lambda}_{2}-\wt{\lambda}_{1}\|_{L^\infty(\ren)}\|u_1\|_{L^2(\dom\cup\bnd)}
    \|u_2-u_1\|_{L^2(\dom)}
    +\|K\|_{L^\infty(\dom)}\|f\|_{L^2(\dom)}\|u_2-u_1\|_{L^2(\dom)}
\end{equation}
or
\begin{equation}\label{E:BndIIb}
    |II|
    \le
    \|\wt{\mu}_2-\wt{\mu}_1\|_{L^2(\ren\times\ren)}\|u_1\|_{L^2(\dom\cup\bnd)}
    \|u_2-u_1\|_{L^2(\dom)}
    +\|K\|_{L^\infty(\dom)}\|f\|_{L^2(\dom)}\|u_2-u_1\|_{L^2(\dom)}.
\end{equation}

Finally, we now combine the bounds for $I$ and $II$ to conclude the proof. Assuming $\wt{M}_2<C_P^{-1}$, we may absorb the bound for $|I|$ given by~\eqref{E:BndIa} into the lower bound in~\eqref{E:MainIneq}. Then using~\eqref{E:BndIIa}, we get
\begin{multline*}
    \|u_2-u_1\|_{L^2(\dom)}^2\le\frac{C_P}{1-C_P\wt{M}_2}
    \ls\|\wt{\gamma}_{1,2}\|_{L^\infty(\ren)}
    \|\wt{\lambda}_{2}-\wt{\lambda}_{1}\|_{L^\infty(\ren)}\|u_1\|_{L^2(\dom\cup\bnd)}
         \|u_2-u_1\|_{L^2(\dom)}\rt\\
    \lt\phantom{\|\wt{\lambda}_{2}-\wt{\lambda}_{1}\|_{L^\infty(\ren)}}
   +\|K\|_{L^\infty(\dom)}\|f\|_{L^2(\dom)}\|u_2-u_1\|_{L^2(\dom)}\rs.
\end{multline*}
The part (a) of the theorem follows upon dividing both sides of the inequality by $\|u_2-u_1\|_{L^2(\dom)}$. The argument for part (b) is similar, using~\eqref{E:BndIb} and~\eqref{E:BndIIb}.
\end{proof}

\begin{rem}
As a consequence of Theorem \ref{thm:kernelBound} we can extract the dependence of nonlocal solutions on different features of the kernel, such as the size of the support of interaction and degree of (integrable) singularity. 
 Let $\mu_1,\mu_2$ be symmetric kernels with support in $B_{\delta_1}(x)$, respectively in $B_{\delta_2}(x)$, and such that  
 \[
 0<m_1< \|\mu_i\|_{L^1(B_{\delta_i})}<m_2<\infty, \quad i=1,2. 
 \]
 Then
\begin{align*}
    \|u_2-u_1\|_{L^2(\dom)} &\leq  C_P
    \ls2\|\wt{\mu}_2-\wt{\mu}_1\|_{L^2(\ren\times\ren)}\|u_1\|_{L^2(\dom\cup\bnd)}
    +
   \|K\|_{L^\infty(\dom)}\|f\|_{L^2(\dom)}\rs\\
   &= C
    \| \,\mu_2\|\mu_1\|_{L^1(B_{\delta_1})}-\mu_1\|\mu_2 \,\|_{L^1(B_{\delta_2})}\|_{L^2(\ren\times\ren)})\|\mu_1\|_{L^1(B_{\delta_1})}^{-1}\|\mu_2\|_{L^1(B_{\delta_2})}^{-1} \\
    &+
   C\left| \|\mu_1\|_{L^1(B_{\delta_1})}-\|\mu_2\|_{L^1(B_{\delta_2})} \right| \|\mu_1\|_{L^1(B_{\delta_1})}^{-1} \|\mu_1\|_{L^1(B_{\delta_1})}^{-1},
\end{align*}
where $C$ depends on the Poincar\'e constant, $u_1$ and $f$. It can be easily shown that if the kernels are simply the characteristic functions of the balls of radii $\delta_1$, respectively $\delta_2$, (i.e. $\mu_1(x,y)=\chi_{B_{\delta_1}}(x-y), \, \mu_2(x,y)=\chi_{B_{\delta_2}}(x-y)$), then one can show that for horizons $\delta_1, \delta_2$ bounded below ($m<\delta_1<\delta_2$, for some $m>0$) we have
\[
 \|u_2-u_1\|_{L^2(\dom)} \leq C(m) |\delta_2-\delta_1|.
\]
The results of \cite[Section 3, Prop. 3.2] {burkovska2020affine} also prescribe that solutions have a Lipschitz variation with respect to the size of $\delta$.
\end{rem}

\section{Continuous Dependence of the Nonlocal Boundary Value Problem in the Nonlinear Setting} \label{nonlinearsec}

We consider two different instances of nonlinearities in terms of continuous dependence. The first is a direct result of Theorem 3.1 and includes nonlinearites inside the nonlocal Laplacian and yields continuous dependence on the forcing term.
\begin{corollary}\label{nonlinear1}
Consider the problem
\begin{equation}
    \left\{\begin{array}{ll}
    \displaystyle{\Int_{\Omega \cup \Gamma}(h(u_i(y))-h(u_i(x)))\mu(x,y)dy=f_i(x)},
        & \quad x \in \Omega,\\
    u_i(x)=g(x),
        & \quad x \in \Gamma.
    \end{array}\right.
\end{equation}
for $i=1,2$, where the nonlinearity $h$ satisfies that $|z_2-z_1| \leq C|h(z_2)-h(z_1)|$.  Let $r \geq 1$. Let $\gamma_\mu,\lambda_\mu\in L^\infty(\ren)$, we have
\[
\|u_1-u_2\|_{L^{r}}< C\|f_1-f_2\|_{L^{r}}. \]
The constant $C$ above is given by
\begin{equation}\label{const3}
C:\frac{\|\frac{1}{\lambda}\|_{L^\infty(\Omega \cup \Gamma)}}{1- M_{\mu,r}\|\frac{1}{\lambda}\|_{L^{\infty}(\Omega \cup \Gamma)}}.
\end{equation}
\end{corollary}

\begin{proof}
Let $h$ satisfy that $|z_2-z_1| \leq |h(z_2)-h(z_1)|$ and define $v_i:=h(u_i)$. Then we are instead considering the problem
\begin{equation}
    \left\{\begin{array}{ll}
    \displaystyle{\Int_{\Omega \cup \Gamma}(v_i(y)-v_i(x))\mu(x,y)dy=f_i(x)},
        & \quad x \in \Omega,\\
    v_i(x)=h(g(x)),
        & \quad x \in \Gamma.
    \end{array}\right.
\end{equation}
for $i=1,2$. Using Theorem \ref{MVPfasym}, we then have
\[\|v_2-v_1\|_{L^{r}(\Omega)}\leq C\|f_2-f_1\|_{L^{r}(\Omega)}.\]

And so, since $\|u_2-u_1\|_{L^{2p}} \leq \|h(u_2)-h(u_1)\|_{L^{2p}}$, we have that
\[\|u_2-u_1\|_{L^{r}(\Omega)}\leq C\|f_2-f_1\|_{L^{r}(\Omega)}.\]

\end{proof}

\begin{corollary}\label{nonlinear2}
Consider the problem
\begin{equation}
    \left\{\begin{array}{ll}
    \displaystyle{\Int_{\Omega \cup \Gamma}(h(u_i(y))-h(u_i(x)))\mu(x,y)dy=f_i(x)},
        & \quad x \in \Omega,\\
    u_i(x)=g(x),
        & \quad x \in \Gamma.
    \end{array}\right.
\end{equation}
for $i=1,2$, where the nonlinearity $h$ satisfies that $|z_2-z_1| \leq C|h(z_2)-h(z_1)|$. Let $g\in L^2(\Gamma)$.
Then, if (M3) and (M4) are satisfied and $1>M_{\mu_{\text{asym}},2}C_P$,
\[ \| u_2-u_1\|_{L^2(\Omega)}
    \leq\frac{C_P}{1-M_{\mu_{\text{asym}},2}C_P}\|f_2-f_1\|_{L^2(\Omega)},
    \]
where $C_P$ is the Poincaré constant from \cite{foss2019nonlocal}.
\end{corollary}

\begin{proof}
Let $h$ satisfy that $|z_2-z_1| \leq |h(z_2)-h(z_1)|$ and define $v_i:=h(u_i)$. Then we are instead considering the problem
\begin{equation}
    \left\{\begin{array}{ll}
    \displaystyle{\Int_{\Omega \cup \Gamma}(v_i(y)-v_i(x))\mu(x,y)dy=f_i(x)},
        & \quad x \in \Omega,\\
    v_i(x)=h(g(x)),
        & \quad x \in \Gamma.
    \end{array}\right.
\end{equation}
for $i=1,2$. Using the above Theorem \ref{linearforcingenergy}, we then have
\[\|v_2-v_1\|_{L^{2}(\Omega)}\leq C_P\|f_2-f_1\|_{L^{2}(\Omega)}.\]

And so, since $\|u_2-u_1\|_{L^{2p}} \leq \|h(u_2)-h(u_1)\|_{L^{2p}}$, we have that
\[\|u_2-u_1\|_{L^{2p}(\Omega)}\leq C_P\|f_2-f_1\|_{L^{2p}(\Omega)}.\]

\end{proof}

\begin{ex}
Consider the following nonlocal boundary value problem:
\begin{equation}
    \left\{\begin{array}{ll}
    \displaystyle{\Int_{(-\delta, 1+\delta)}(\sin(u(y))-\sin(u(x)))\mu(x,y)dy=f_i(x)},
        & \quad x \in \Omega ,\\
    u_i(x)=g(x),
        & \quad x \in \Gamma.
    \end{array}\right.
\end{equation}
From Corollaries \ref{nonlinear1} and \ref{nonlinear2}, we  know that 

\[
\|u_1-u_2\|_{L^{r}}< C\|f_1-f_2\|_{L^{r}}. \]
The constant $C$ above is given by
\begin{equation*}
C:=\frac{\|\frac{1}{\lambda}\|_{L^\infty(\Omega \cup \Gamma)}}{1- M_{\mu,r}\|\frac{1}{\lambda}\|_{L^{\infty}(\Omega \cup \Gamma)}}
\end{equation*}
or if (M3) and (M4) are satisfied and $1>M_{\mu_{\text{asym}},2}C_P$,
\[ \| u_2-u_1\|_{L^2(\Omega)}
    \leq\frac{C_P}{1-M_{\mu_{\text{asym}},2}C_P}\|f_2-f_1\|_{L^2(\Omega)}\]
where $C_P$ is the Poincaré constant from \cite{foss2019nonlocal}.
\end{ex}

In the following we will study continuous dependence of solutions on the profiles of nonlinearities appearing in the forcing term. 
\begin{theorem}\label{nonlinear3}
Consider the nonlocal systems
\begin{equation}
    \left\{\begin{array}{ll}
    \displaystyle{\mathcal{L}_{\mu}u_i(x)=f_i(x,u_i(x))},
        & \quad x \in \Omega,\\
    u_i(x)=g(x),
        & \quad x \in \Gamma.
    \end{array}\right.
\end{equation}
for $i=1,2$, where $f_1, f_2$ are Lipschitz in $u$, i.e. there exist $L_1, L_2>0$ such that 
\[
|f_i(x,u)-f_i(x,v)|\leq L_i(x)|u-v|, \quad x, u, v\in \re.
\]
 Then if (M3) and (M4) are satisfied and $1-M_{\mu_{\text{asym}},2}C_P-L_1C_P>0$, where $C_P$ is the Poincar\'e constant and where $L_1(x)$ is the Lipschitz constant associated with $f_1$ with respect to the first argument and $L_1=\|L_1\|_{L^1(\Omega)}$. 
\[\| u_2-u_1\|_{L^2(\Omega)} \leq  C\|f_2-f_1\|_{L^\infty(\Omega \times\mathbb{R})}\]
where $C=\frac{C_P}{1-M_{\mu_{\text{asym}},2}C_P-L_1C_P}$.
\end{theorem}

\begin{proof}
From Theorem 3.5, we have that 
\begin{align*}
   \| u_2-u_1\|_{L^2(\Omega)}
    &\leq\frac{C_P}{1-M_{\mu_{\text{asym}},2}C_P}\|f_2(\cdot,u_2(\cdot))-f_1(\cdot,u_1(\cdot))\|_{L^2(\Omega)}.
\end{align*}
Using the boundedness of $f_1,f_2$ we have
\begin{align*}
    \|f_2(\cdot,u_2(\cdot))-f_1(\cdot,u_1(\cdot))\|_{L^2(\Omega)} &\leq\|f_2(\cdot,u_2(\cdot))-f_1(\cdot,u_2(\cdot))\|_{L^2(\Omega)}+\|f_1(\cdot,u_2(\cdot))-f_1(\cdot,u_1(\cdot))\|_{L^2(\Omega)}\\
    &\leq\|f_2-f_1\|_{L^\infty(\Omega,\mathbb{R})}+\|f_1(\cdot,u_2(\cdot))-f_1(\cdot,u_1(\cdot))\|_{L^2(\Omega)}.
\end{align*}
Since $f_1$ is Lipschitz in $u$, we have
\begin{align*}
  \|f_2(\cdot,u_2(\cdot))-f_1(\cdot,u_1(\cdot))\|_{L^2(\Omega)}  &\leq\|f_2-f_1\|_{L^\infty(\Omega,\mathbb{R})}+\|L_1(\cdot)(u_2(\cdot)-u_1(\cdot))\|_{L^2(\Omega)}\\
   &\leq\|f_2-f_1\|_{L^\infty(\Omega,\mathbb{R})}+L_1\|u_2-u_1\|_{L^2(\Omega)}.
\end{align*}
Hence
\begin{align*}
    \| u_2-u_1\|_{L^2(\Omega)} \leq  \frac{C_P}{1-M_{\mu_{\text{asym}},2}C_P}\Big( \|f_2-f_1\|_{L^\infty(\Omega,\mathbb{R})}+L_1\|u_2-u_1\|_{L^2(\Omega)} \Big),
\end{align*}
or
\begin{align*}
    \| u_2-u_1\|_{L^2(\Omega)} \leq \frac{C_P}{1-M_{\mu_{\text{asym}},2}C_P-L_1C_P} \|f_2-f_1\|_{L^\infty(\Omega,\mathbb{R})}.
\end{align*}
Since $1-M_{\mu_{\text{asym}},2}C_P-L_1C_P>0$, the conclusion giving sensitivity of solutions with respect to forcing follows.
\end{proof}


\section{Numerical Results} \label{sec:Numerical}
In this section we numerically solve 
\begin{subequations}
  \label{eq:nonlocalLaplacian}
  \begin{empheq}[]{alignat=2}
    \int_{x-\delta}^{x+\delta} \left( u(y) - u(x) \right) \mu(x,y) \, dy  &= f(x) \quad\quad & x\in \Omega \\
    u(x) &= g(x) \quad\quad &x\in \Gamma
  \end{empheq}
\end{subequations}
to illustrate the bounds presented thus far. The domain $\Omega$ is chosen to be $(0,1)$, so that $\Omega \cup \Gamma$ for all examples is $(-\delta,1+\delta)$. All numerical results are computed using the discontinuous Galerkin discretizaton described in  \cite{chen2011continuous} and using a uniform mesh spacing of $h=1/200$.

\subsection{Sensitivity with respect to perturbations in the forcing term}
We explore how perturbations to the the right-hand side $f(x)$ perturb the solution. Given the ability of nonlocal methods to handle discontinuous solutions, we consider two separate forcing functions that will produce discontinuous solutions. From Theorem \ref{linearforcingenergy}, we know that if $1>M_{\mu_{\text{asym}},2}\,C_P$, then
\[ \| u_2-u_1\|_{L^2(\Omega)}
    \leq\frac{C_P}{1-M_{\mu_{\text{asym}},2}C_P}\|f_2-f_1\|_{L^2(\Omega)},\]
where $C_P$ is the Poincaré constant from Lemma \ref{Poinc}. If we select a symmetric kernel, then $M_{\mu_{\text{asym}},2}=0$ and the above inequality reduces to 
\begin{equation} \label{eq:forcingBound}
\| u_2-u_1\|_{L^2(\Omega)} \leq C_P\|f_2-f_1\|_{L^2(\Omega)}.
\end{equation}
For both of the following cases, as in Remark \ref{rem:poincareConstant}, we let $\mu_\delta\equiv3\delta^{-3}$ on $B_\delta(x)$, so $C_P=\frac{9}{8}$ for all $\delta>0$. We will compute the norms in \eqref{eq:forcingBound} and verify that it is satisfied numerically for these specific examples.

\subsubsection{Discontinuous forcing sinusoidal perturbation} \label{sec:sinusoid}
\begin{figure}
    \centering
    \begin{subfigure}[b]{0.45\textwidth}
      \centering
      \includegraphics[scale=0.25]{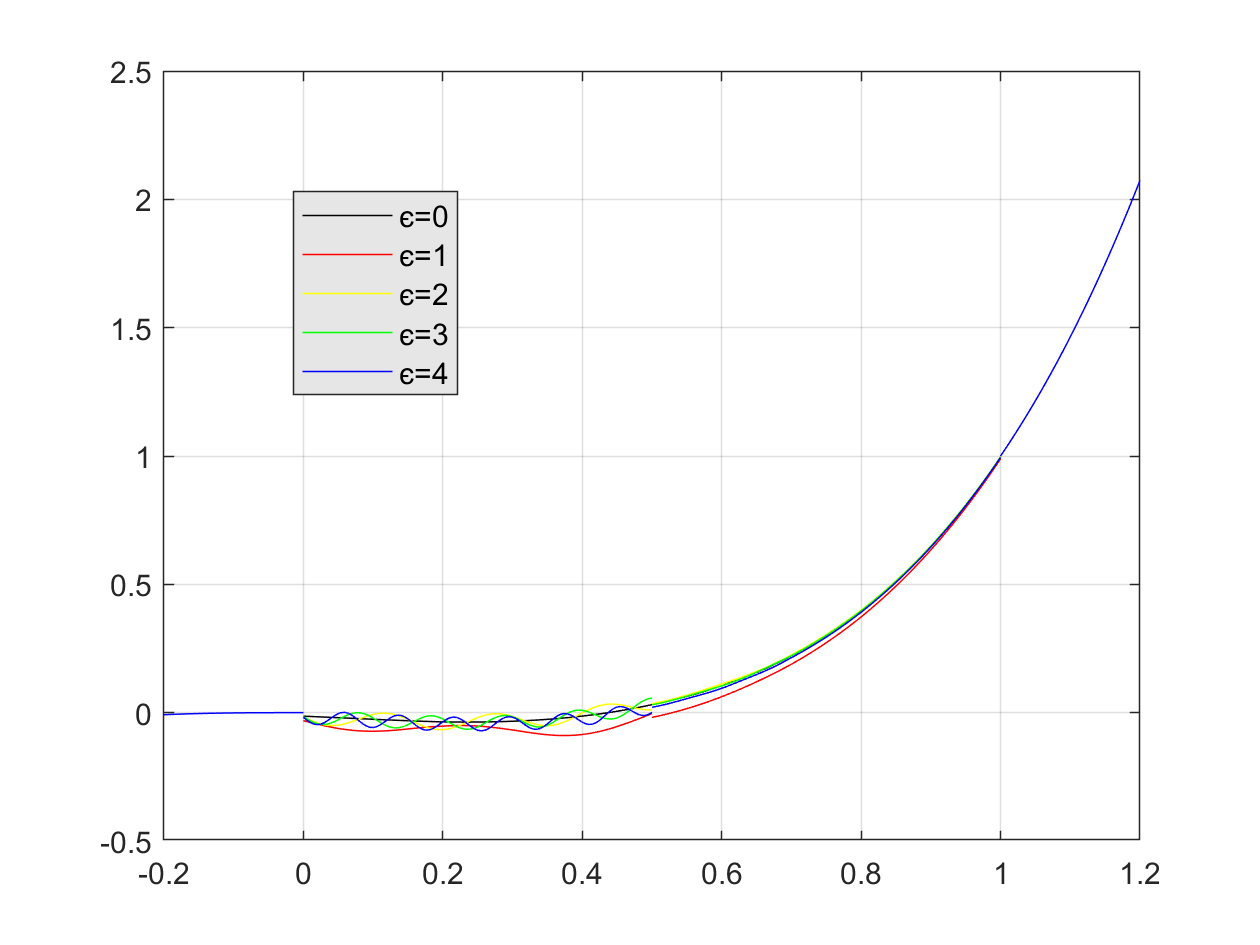}
      \caption{}
      \label{fig:sinForcing}
    \end{subfigure}
    \hfill
    \begin{subfigure}[b]{0.45\textwidth}
      \centering
      \includegraphics[scale=0.25]{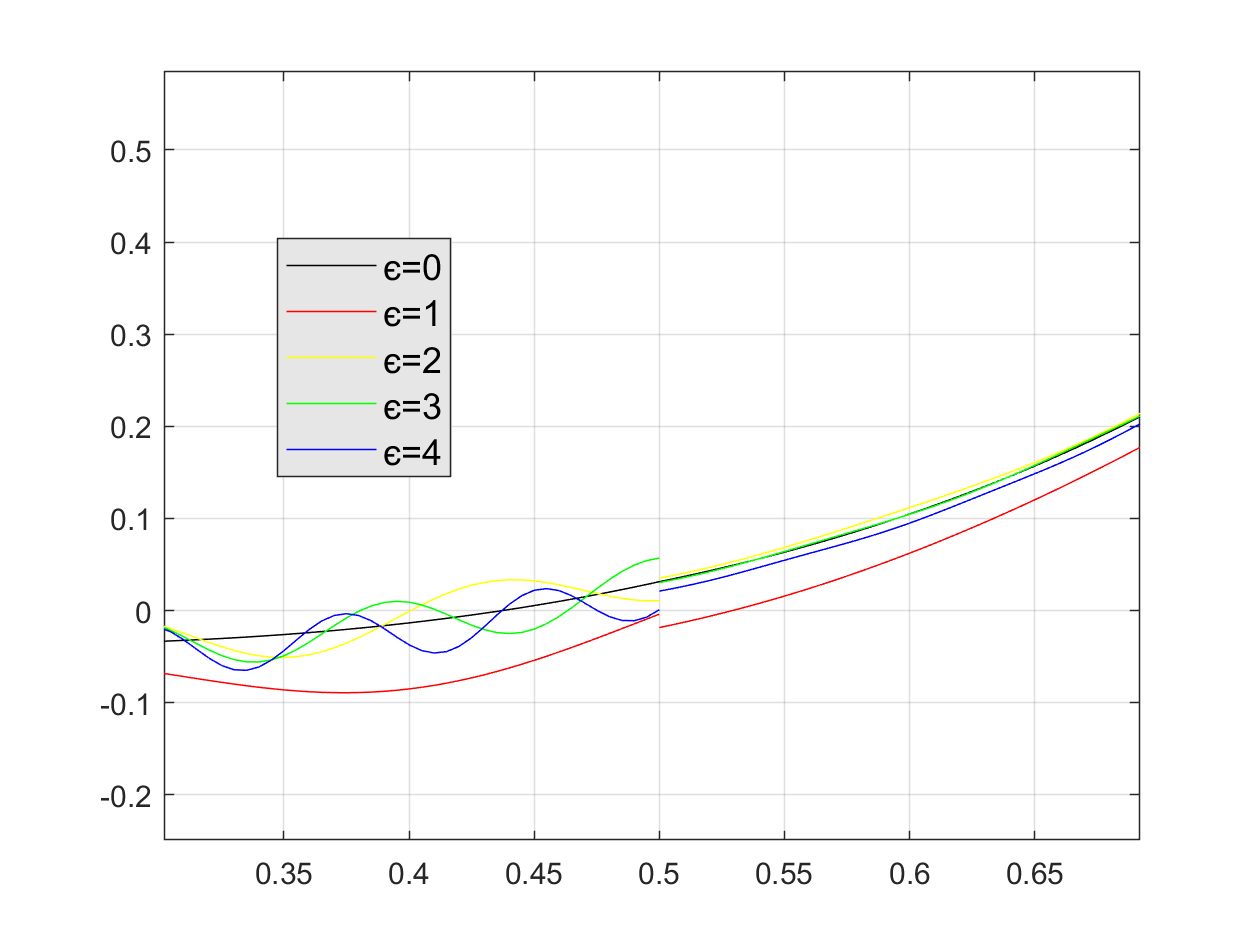}
      \caption{}
    \end{subfigure}
    \par \bigskip
    \begin{subtable}{0.45\textwidth}
      \centering
      \begin{tabular}{|c|c|c|c|}
        \hline
        $\varepsilon$	&$\|f_{\varepsilon}-f_0 \|_{L^2}$& $\|u_{\varepsilon}-u_0 \|_{L^2}$ &	 $\frac{\|u_{\varepsilon}-u_0 \|_{L^2}}{\|f_{\varepsilon}-f_0 \|_{L^2}}$ \\
        \hline
        1.0 &1.95421 &0.0385269 &0.019715 \\
        \hline
        2.0	&1.98031 & 0.0150158&0.0075826 \\
        \hline
        3.0	&1.99662 & 0.013236&0.0066292 \\
        \hline
        4.0 &1.98945 &0.015445 &0.0077636 \\
        \hline
      \end{tabular}
    \caption{} 
    \label{table:sinusoid}
    \end{subtable}
    \caption{Numerical results for the example of \S\ref{sec:sinusoid}. (a): Solutions $u_\varepsilon(x)$ over $(-\delta,1+\delta)$ for various $\varepsilon$ using the piecewise sinusoidal forcing function \eqref{eq:sineforcing}. (b): A zoom in to show the discontinuity at $x=0.5$. (c): Numerical evaluation of the terms in \eqref{eq:forcingBound} for various $\varepsilon$. For these particular examples the ratio in the rightmost column in the table is less than 9/8, the Poincar\'{e} constant in \eqref{eq:forcingBound}, showing that \eqref{eq:forcingBound} is satisfied.}
\end{figure}

In this example $\delta=0.2$, the collar conditions are $g(x)=x^3$ in $(-\delta,0)$ and $g(x)=x^4 $ in $(1,1+\delta)$, and the forcing function is
\begin{equation}\label{eq:sineforcing}
f_{\varepsilon}(x)=\begin{cases}
6x+ 4\sin(20\varepsilon x) &x\leq 0.5\\
12x^2 & x>0.5.
\end{cases}    
\end{equation}
We study the sensitivity of the solutions on the forcing term by varying the parameter $\varepsilon$. Note that $f_{\varepsilon=0}(x)$ is continuous, but $f_\varepsilon(x)$ is not continuous for $\varepsilon>0$. In this example a discontinuity in the forcing function is sufficient to force a discontinuity in the solution. In Figure~\ref{fig:sinForcing}, we plot $u_{\varepsilon}(x)$ for various $\varepsilon$, and in Table \ref{table:sinusoid} we show numerically that \eqref{eq:forcingBound} is satisfied for this example.

\subsubsection{Sigmoid Forcing}\label{sec:sigmoid}

\begin{figure}
    \centering
    \begin{subfigure}[b]{0.45\textwidth}
      \centering
      \includegraphics[scale=0.25]{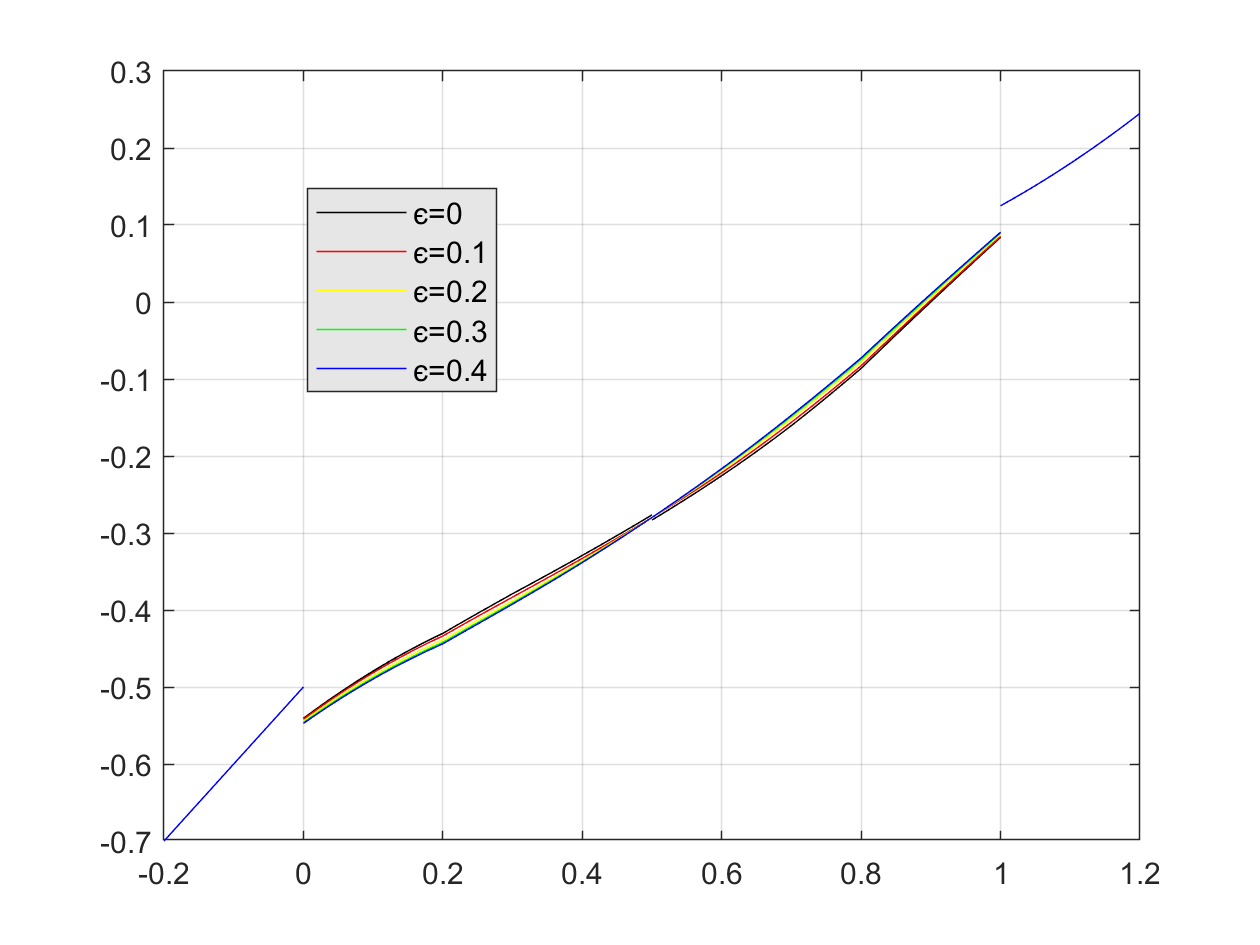}
      \caption{}
      \label{fig:sigmoid}      
    \end{subfigure}
    \hfill
    \begin{subfigure}[b]{0.45\textwidth}
      \centering
      \includegraphics[scale=0.25]{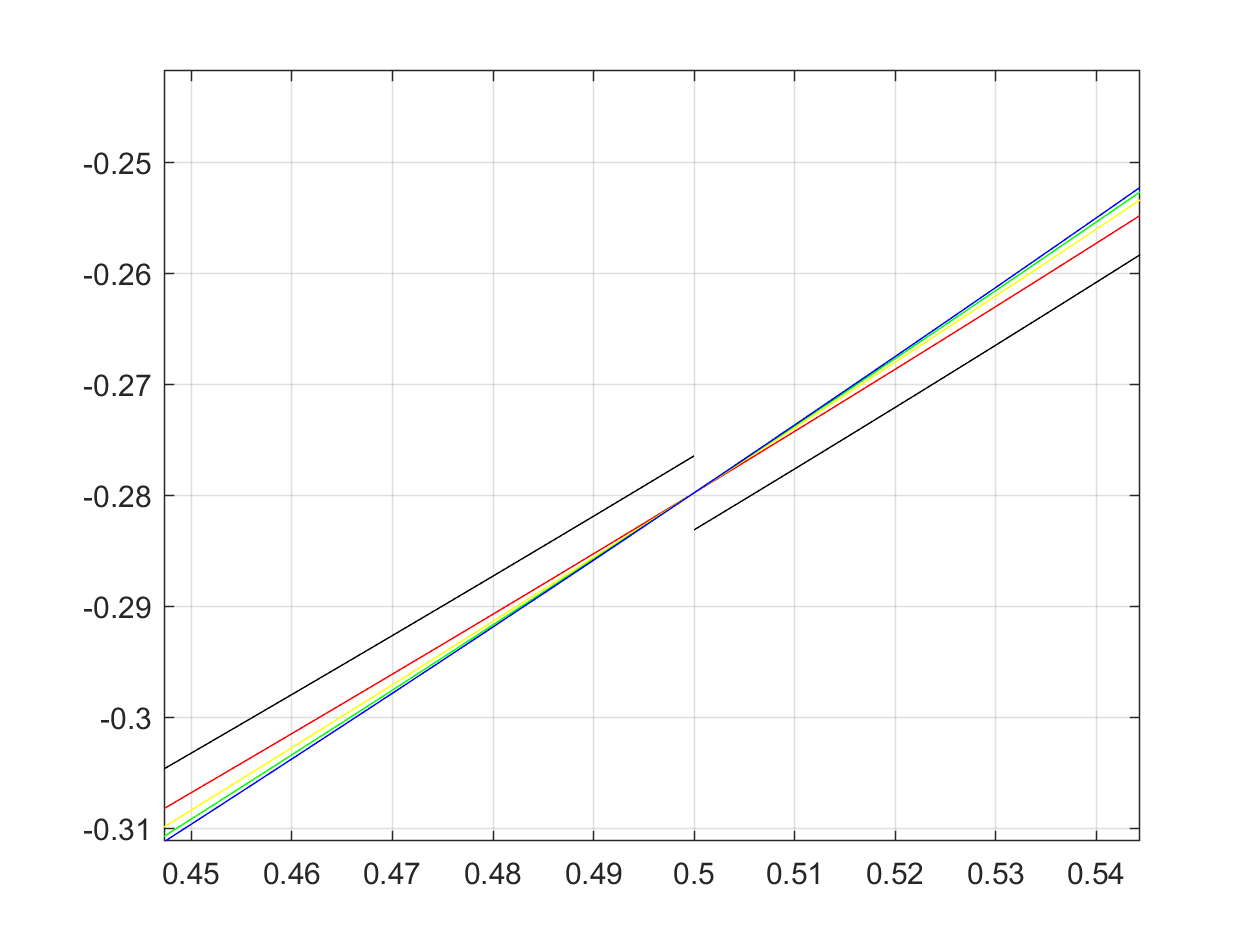}
      \caption{}
    \end{subfigure}
    \par \bigskip
    \begin{subtable}{0.45\textwidth}
      \centering
      \begin{tabular}{|c|c|c|c|}
        \hline
        $\varepsilon$	&$\|f_{\varepsilon}-f_0 \|_{L^2}$& $\|u_{\varepsilon}-u_0 \|_{L^2}$ &	 $\frac{\|u_{\varepsilon}-u_0 \|_{L^2}}{\|f_{\varepsilon}-f_0 \|_{L^2}}$ \\
\hline
0.1 &0.194463 &0.00333435 &0.0171464 \\
\hline
0.3 &0.274272 &0.00705121 &0.0257089 \\
\hline
0.3 & 0.326465&0.00937336 &0.0287117 \\
\hline
0.4 &0.360981 &0.0108048 & 0.0299318\\
        \hline
      \end{tabular}
    \caption{} 
    \label{table:sigmoid}
    \end{subtable}
    \caption{Numerical results for the example of \S\ref{sec:sigmoid}. (a): Solutions $u_\varepsilon(x)$ over $(-\delta,1+\delta)$ for various $\varepsilon$ using the sigmoid forcing function \eqref{eq:sigmoid}. (b): A zoom in to show the discontinuity at $x=0.5$. (c): Numerical evaluation of the terms in \eqref{eq:forcingBound} for various $\varepsilon$. For these particular examples the ratio in the rightmost column in the table is less than 9/8, the Poincar\'{e} constant in \eqref{eq:forcingBound}, showing that \eqref{eq:forcingBound} is satisfied.}
\end{figure}

In this example $\delta=0.2$, the collar conditions are $g(x)=x-0.5$ in $(-\delta,0)$ and $g(x)=\frac{(x-0.5)^2}{2}$ in $(1,1+\delta)$, and the forcing function is
\begin{equation}\label{eq:sigmoid}
f_{\varepsilon}(x)=\frac{e^{\frac{x-0.5}{\varepsilon}}}{1+e^{\frac{x-0.5}{\varepsilon}}}.
\end{equation}
We study the sensitivity of the solutions on the forcing term by varying the parameter $\varepsilon$. Note that $f_{\varepsilon=0}(x)$ is discontinuous, but $f_\varepsilon(x)$ is continuous for $\varepsilon>0$. In Figure~\ref{fig:sigmoid}, we plot $u_{\varepsilon}(x)$ for various $\varepsilon$, and in Table \ref{table:sigmoid} we show numerically that \eqref{eq:forcingBound} is satisfied for this example. 

\subsection{Sensitivity with respect to perturbations in the boundary data}\label{sec:collar}

From Theorem \ref{collarenergy}, we know that if $1>M_{\mu_{\text{asym}},2}\,C_P$, then 
\begin{equation}
    \| u_2-u_1\|_{L^2(\Omega)}\leq 
   \frac{C_P\|\mu\|_{L^2(\Omega\times \Gamma)}}{ 1-C_PM_{\mu_{\text{asym}},2}}\|g_2-g_1\|_{L^2(\Gamma)},
\end{equation}
where $C_P$ is the Poincaré constant from \cite{foss2019nonlocal}. If we select symmetric kernels, $M_{\mu_{\text{asym}},2}=0$, then the above inequality reduces to 
\begin{equation} \label{eq:collarBound}
\| u_2-u_1\|_{L^2(\Omega)} \leq C_P\|\mu\|_{L^2(\Omega\times \Gamma)}\|g_2-g_1\|_{L^2(\Gamma)}.
\end{equation}
Let $\mu(x,y) \equiv 3\delta^{-3}$ on $B_\delta(x)$ so that  $C_P=\frac{9}{8}$ for all $\delta>0$ as in Remark \ref{rem:poincareConstant}. We observe that $\|\mu\|_{L^2(\Omega\times \Gamma)}=6\delta^{-2}$. We chose the horizon $\delta=0.1$, which means $C_P\|\mu\|_{L^2}$, our constant of proportionality, becomes $675$.

\begin{figure}
    \centering
    \begin{subfigure}[c]{0.45\textwidth}
      \centering
      \includegraphics[scale=0.25]{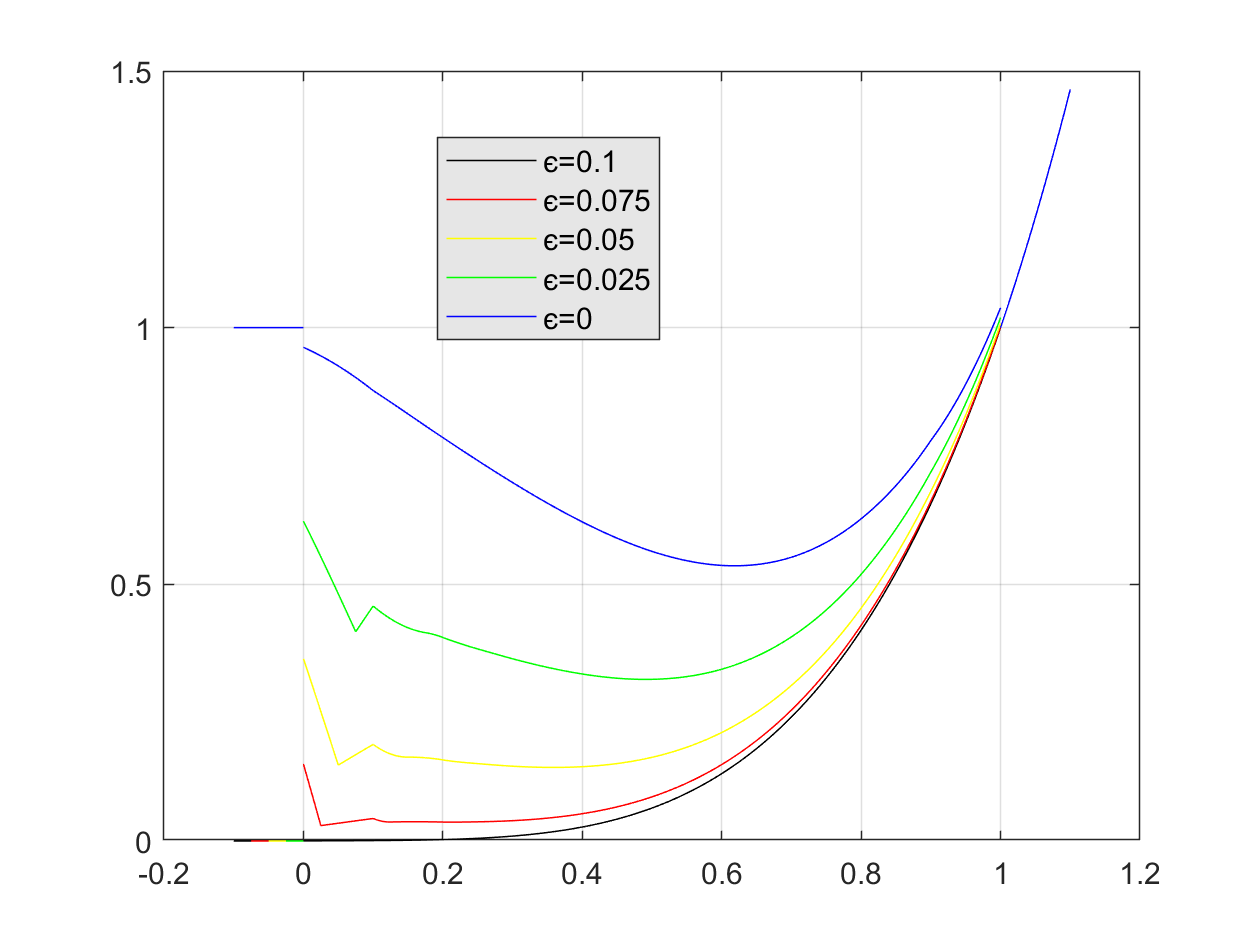}
      \caption{}
      \label{fig:discontcollar}
    \end{subfigure}
    \hfill
    \begin{subtable}[c]{0.45\textwidth}
      \centering
      \begin{tabular}{|c|c|c|c|}
        \hline
        $\varepsilon$	&$\|g_{\varepsilon}-g_{\delta} \|_{L^2}$& $\|u_{\varepsilon}-u_{\delta} \|_{L^2}$ &	 $\frac{\|u_{\varepsilon}-u_{\delta} \|_{L^2}}{\|g_{\varepsilon}-g_{\delta} \|_{L^2}}$ \\
        \hline
        0.075&0.158104 &0.160571 &1.0156 \\
        \hline
        0.05&0.223598 &0.253197 &1.13237 \\
        \hline
        0.025&0.273854&0.400205 &1.46138 \\
        \hline
        0&0.316221 &0.65106 &  2.05887\\
        \hline
      \end{tabular}
    \caption{} 
    \label{table:discontcollar}
    \end{subtable}
    \caption{Numerical results for the example of \S\ref{sec:collar}. (a): Solutions $u_\varepsilon(x)$ over $(-\delta,1+\delta)$ for various $\varepsilon$ using the collar given by \eqref{eq:piecewisecollar}. (b): Numerical evaluation of the terms in \eqref{eq:collarBound} for various $\varepsilon$. For these particular examples the ratio in the rightmost column in the table is less than 675, showing that \eqref{eq:collarBound} is satisfied.}
\end{figure}

With the forcing given by $f(x)=12x^2$, we vary parameter $\varepsilon>0$ and consider the collar data
\begin{equation}\label{eq:piecewisecollar}
g_\varepsilon(x)=\begin{cases}1, & x\in (-\delta, -\varepsilon) \\
x^4, & x\in  [-\varepsilon,0) \cup (1, 1+\delta) 
\end{cases}.
\end{equation}
If $\varepsilon = \delta$ or $\varepsilon = 0$ then $g_{\varepsilon=\delta}(x)$ and $g_{\varepsilon =0}(x)$ are continuous, but otherwise $g_\varepsilon(x)$ is discontinuous. We illustrate the sensitivity of the solutions on the boundary data by varying the parameter $\varepsilon$. In Figure~\ref{fig:discontcollar}, we plot $u_{\varepsilon}(x)$ for various $\varepsilon$, and in Table \ref{table:discontcollar} we show numerically that \eqref{eq:collarBound} is satisfied for this example. 

\subsection{Sensitivity with respect to perturbations in the kernel}
Recalling Theorem \ref{thm:kernelBound}, it was proven that the analytical bound for perturbations in the kernel is 
\begin{equation}
\label{eq:generalkernelBound} \|u_2-u_1\|^2_{L^2(\dom)}
    \le
    \frac{C_P}{1-C_P\|\wt{\mu}_{2,\asym}\|_{L^2(\dom\times\dom)}}
    \left(\|u_1\|_{L^2(\dom)}
    \|\wt{\mu}_2-\wt{\mu}_1\|_{L^2(\dom\times\dom)}
    +\|K\|_{L^\infty(\re^n)}\|f\|_{L^2(\dom)}\right),
    \end{equation}
where $M=\left|\frac{1}{\|\mu_2\|_{L^1}}-\frac{1}{\|\mu_1\|_{L^1}}\right|$ and $\wt{\mu}_i=\frac{\mu_i}{\|\mu_i\|_{L^1(\Omega \cup \Gamma)}}$. Note that if we allow only symmetric kernels, our bound becomes 
\begin{equation}
\label{eq:kernelBound}
\|u_2-u_1\|^2_{L^2(\dom)}
    \le
    C_P
    \left(\|u_1\|_{L^2(\dom)}
    \|\wt{\mu}_2-\wt{\mu}_1\|_{L^2(\dom\times\dom)}
    +\|K\|_{L^\infty(\re^n)}\|f\|_{L^2(\dom)}\right).
\end{equation}

\subsubsection{Sensitivity with respect to singularity in the kernel} \label{sec:Kernel}
We study the sensitivity of the solutions upon the the kernel $\mu_\varepsilon(x,y)=\frac{3-\varepsilon}{\delta^{3-\varepsilon}}|x-y|^{-\varepsilon}$ as we vary $\varepsilon\geq0$. We choose the horizon $\delta=0.2>0$, let the forcing function be given by $f(x)=12x^2$, and the collar conditions be given by $g(x)=x^4$. We solve numerically for several values of the parameter $\varepsilon$ and compare them as perturbation to the solution computed using the constant kernel (i.e., $\varepsilon=0$). The numerical solutions for this case are unremarkable, so we show only the tabular data in Table \ref{table:Kernel}, demonstrating that \eqref{eq:kernelBound} is satisfied for this example. Here $C_P=2^{-1}(2-\varepsilon)^{\varepsilon-2}(3-\varepsilon)^{2-\varepsilon}$.

\begin{table}[ht!]
\centering
\begin{tabular}{|c|c|c|c|c|}
   \hline
   $\varepsilon$&$C_P$&$B$ & $\|u_\varepsilon-u_0 \|_{L^2}$ &	 $\frac{\|u_\varepsilon-u_0 \|_{L^2}}{B}$ \\
     \hline
0.2& 1.1076&0.117579 &0.000215331 &0.00120516 \\
\hline
0.4& 1.0873&0.300229 &  0.000447169&0.000987976 \\
\hline
0.6& 1.0634& 0.618071& 0.000697213&0.000756179 \\
\hline
0.8& 1.0348&1.00269 &0.000967517 &0.000651402 \\
\hline
\end{tabular}
\caption{Numerical results for the example of \S\ref{sec:Kernel}. Numerical evaluation of the terms in \eqref{eq:kernelBound} for various $\varepsilon$, where we have defined $B:=  2\|u_0\|_{L^2(\dom)} \|\wt{\mu}_2-\wt{\mu}_1\|_{L^2(\dom\times\dom)} +\|K\|_{L^\infty(\re^n)}\|f\|_{L^2(\dom)}.$ For these particular examples the ratio in the rightmost column in the table is less than the Poincaré constant given in the second column.}
\label{table:Kernel}
\end{table}

\subsubsection{Heterogeneous kernel in $x$}\label{sec:heterogeneousx}
\begin{figure}
    \centering
    \begin{subfigure}[c]{0.6\textwidth}
      \centering
      \includegraphics[scale=0.25]{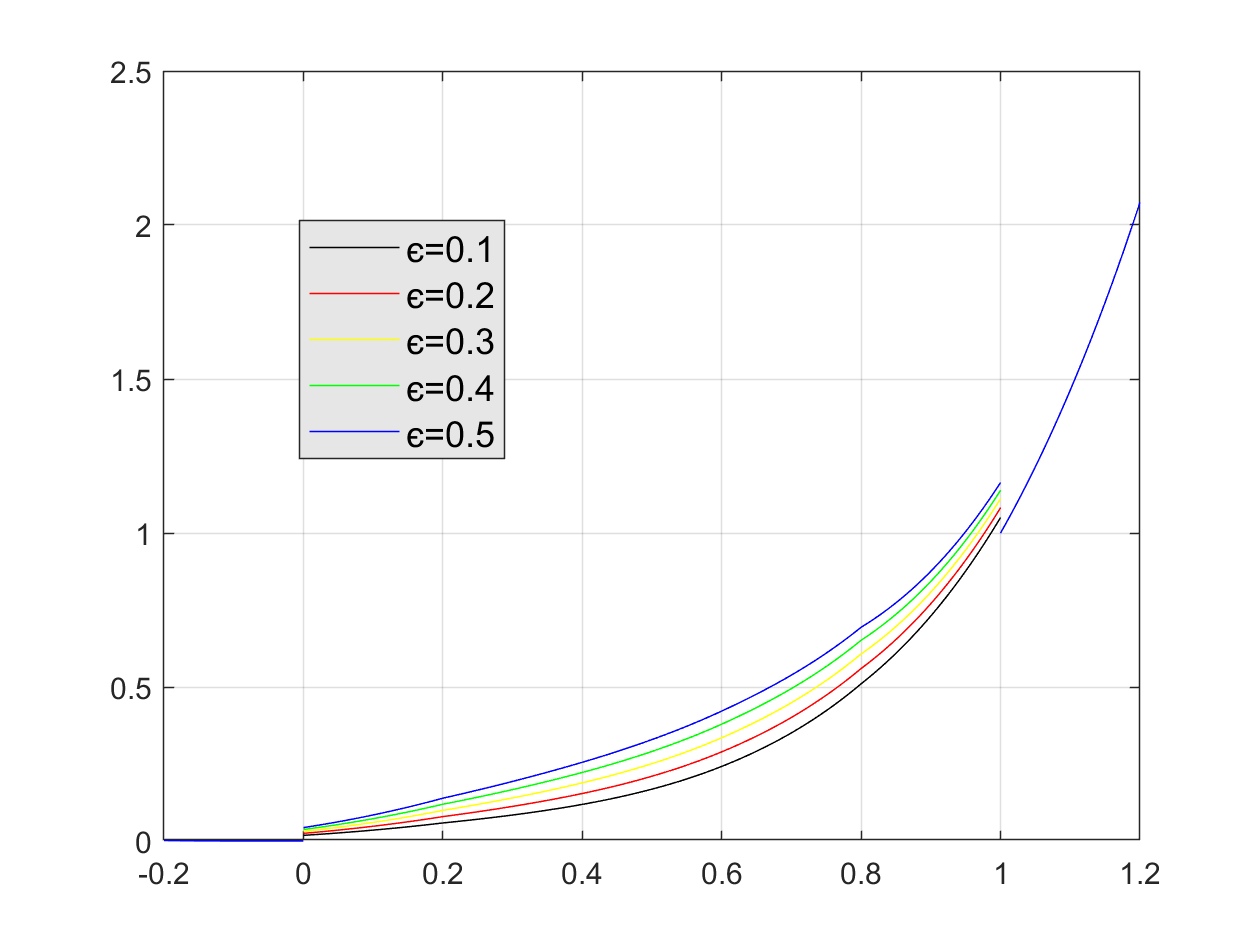}
      \caption{}
      \label{fig:heterogeneousx}
    \end{subfigure}
    \par \bigskip
    \begin{subtable}[c]{0.8\textwidth}
      \centering
      \begin{tabular}{|c|c|c|c|c|c|c|}
      \hline
         $\varepsilon$&$\frac{C_P}{1-C_P\|\wt{\mu}_{\varepsilon,\asym}\|_{L^2(\ren\times\ren)}}$&$\|\widetilde{\mu_\varepsilon}-\widetilde{\mu_0}\|_{L^1(\re^n)}$&$K\|f\|_{L^2(\Omega)}$& $\|u_\varepsilon-u_0 \|_{L^2}$ &	 $\frac{\|u_\varepsilon-u_0 \|_{L^2}}{B}$ \\
\hline
0.2&0.102950&0.0105405 & 93.7685 &0.0369836&0.00039434\\
\hline
0.3&0.0908265&0.0210802& 197.5639&0.0724964& 0.00036694\\
\hline
0.4&0.0801798&0.0316185& 312.4746 &0.106577&0.00034106\\
\hline
0.5&0.0708216&0.0421549 &439.7065 &0.139266&0.00031671\\
\hline
      \end{tabular}
      \caption{}
      \label{table:heterogeneousx}      
    \end{subtable}
   \caption{Numerical results for the example of \S\ref{sec:heterogeneousx}.  (a): Solutions $u_\varepsilon(x)$ over $(-\delta,1+\delta)$ for various $\varepsilon$ found by varying $\varepsilon$ in \eqref{eq:heterogeneouskernel}. (b)Numerical evaluation of the terms in \eqref{eq:generalkernelBound} for various $\varepsilon$, where we have defined $B:=  \|u_1\|_{L^2(\dom)} \|\wt{\mu}_2-\wt{\mu}_1\|_{L^2(\dom\times\dom)} +\|K\|_{L^\infty(\re^n)}\|f\|_{L^2(\dom)}.$ For these particular examples the ratio in the rightmost column is less than the ratio in the second column, demonstrating that \eqref{eq:generalkernelBound} is satisfied for this example.}
\end{figure}

Next we consider a heterogeneous kernel. For varying parameter $\varepsilon\geq0$ and given horizon $\delta=0.2$, let the forcing be given by $f(x)=12x^2$, and the collar data be given by $g(x)=x^4$. We consider the kernel

\begin{equation}\label{eq:heterogeneouskernel}
    \mu(x,y):=
\begin{cases}
\frac{1}{\delta^3}(4-x)e^{xy\varepsilon} \quad |y-x|\leq \delta  \\
0, \quad \text{ otherwise}.
\end{cases}
\end{equation}
As this kernel is nonsymmetric, we must compute the terms in \eqref{eq:generalkernelBound}. We can compute analytically that $\|f\|_{L^2(\Omega)}\approx5.3666$. We define the unperturbed solution $u_1$ as the solution computed when $\varepsilon=0.1$, and compute numerically that $\|u_1\|_{L^2(\Omega \cup \Gamma)}\approx 0.38937$. In Figure~\ref{fig:heterogeneousx}, we plot $u_{\varepsilon}(x)$ for various $\varepsilon$, and in Table \ref{table:heterogeneousx} we show numerically that \eqref{eq:generalkernelBound} is satisfied for this example. 

\subsubsection{Spatially Discontinuous Domain} \label{sec:Spatialdiscont}
We consider a problem where the material in the region $\Xi := (0.5-\varepsilon,0.5+\varepsilon)$ has been removed from $\Omega$, but where $2\varepsilon < \delta = 0.2$, so that the remaining material is still self-connected. This can be realized by removing all bonds between $\Omega$ and $\Xi$, which manifests as a spatially heterogeneous kernel. We study the sensitivity of the solution based on the size of the region removed. We utilize the kernel $\mu(x,y)=3\delta^{-3}$, except in the noted regions where bonds have been removed. We choose a forcing function $f(x)=0$ so that $\|f\|_{L^2(\Omega)}=0$ and let $g(x)=x.$ When $\varepsilon=0$, we have continuity and $\|u_0\|_{L^2(\Omega)}=0.57735$.

\begin{figure}
    \centering
    \begin{subfigure}[b]{0.45\textwidth}
      \centering
      \includegraphics[scale=0.25]{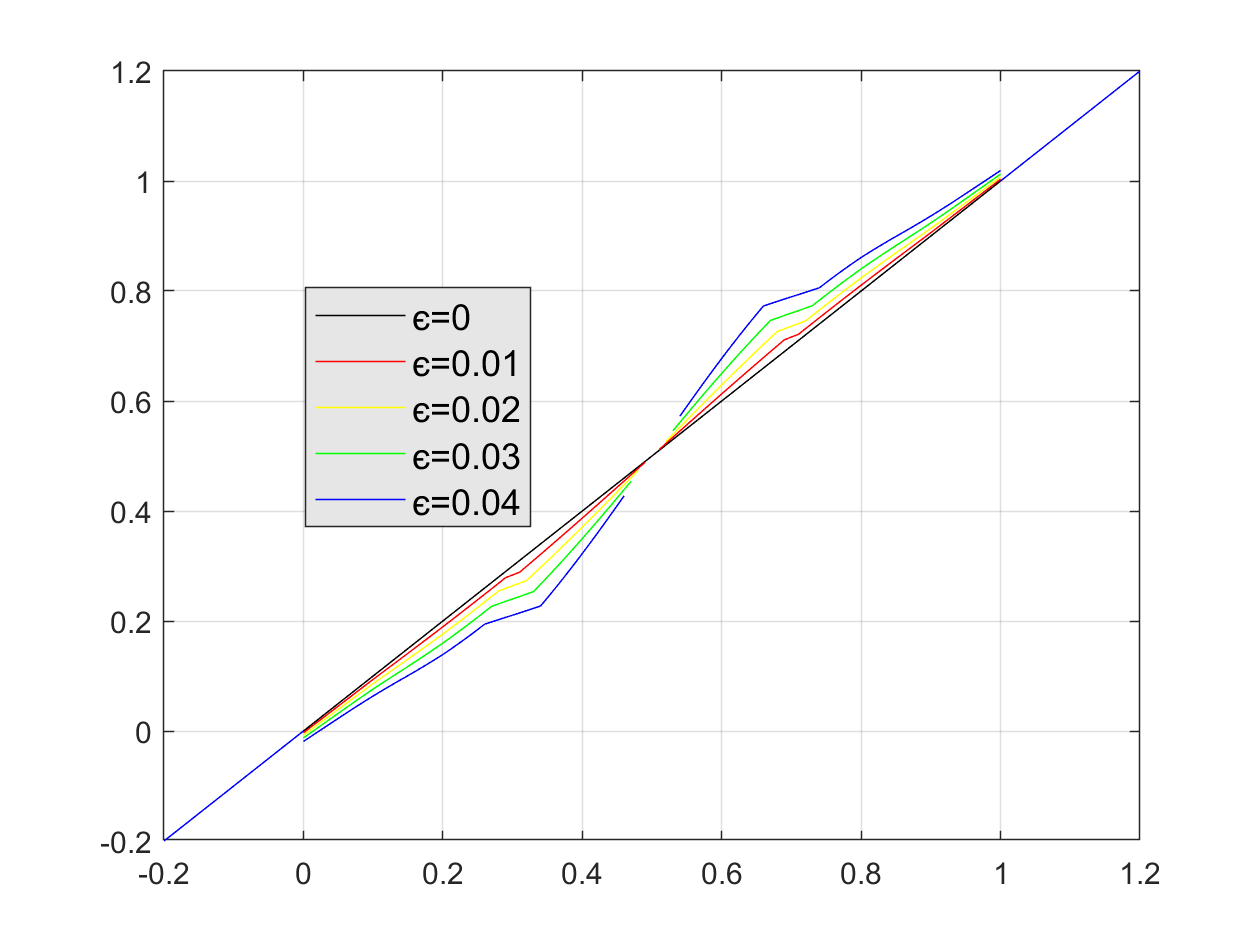}
      \caption{}
      \label{fig:horizon}      
    \end{subfigure}
    \hfill
    \begin{subfigure}[b]{0.45\textwidth}
      \centering
      \includegraphics[scale=0.25]{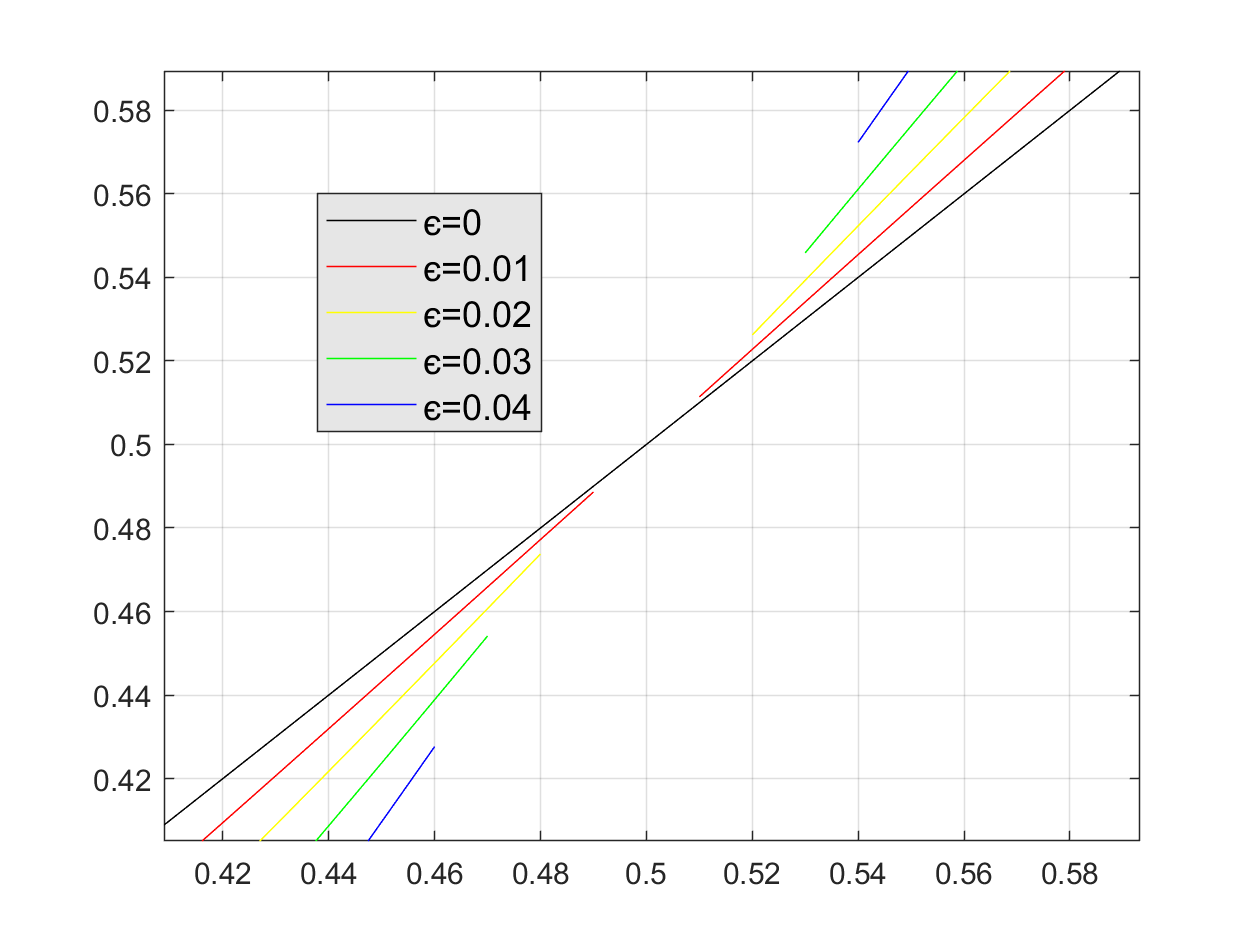}
      \caption{}
    \end{subfigure}
    \par \bigskip
    \begin{subtable}{\textwidth}
      \centering
      \begin{tabular}{|c|c|c|c|c|}
        \hline
        $\varepsilon$
        &
        $\|\widetilde{\mu}_\varepsilon-\widetilde{\mu}_{0}\|_{L^2(\re^n)}$
        &
        $ \|u_2\|_{L^2(\dom\cup\bnd)}$
        &
        $\|u_\varepsilon-u_{0} \|_{L^2}$ 
        &
        $\frac{\|u_\varepsilon-u_{0} \|_{L^2}^2}{B}$ \\
        \hline
        0.01&0.125 & 0.576929 &  0.0715597 &0.03550 \\
        \hline
        0.02&0.25 &0.577764  &  0.102999 &0.0367 \\
        \hline
        0.03&0.375 &0.580346  & 0.129405 &0.0385 \\
        \hline
        0.04&0.5 & 0.585365 &0.154908 & 0.0410\\
        \hline
      \end{tabular}
    \caption{} 
    \label{table:horizon}
    \end{subtable}
    \caption{Numerical results for the example of \S\ref{sec:Spatialdiscont}. (a): Solutions $u_\delta(x)$ over $(-\delta,1+\delta)$ for various $\varepsilon$ found by varying the size of the discontinuity in \eqref{eq:nonlocalLaplacian}. (b): A zoom in to show the differences between the solutions. (c): Numerical evaluation of the terms in \eqref{eq:spatialDiscontBound} for various $\delta$, where we have defined $B:=  2 \|u_2\|_{L^2(\dom\cup\bnd)}
    \|\wt{\mu}_2-\wt{\mu}_1\|_{L^2(\ren\times\ren)}.$
    For these particular examples the ratio in the rightmost column in the table is less than 9/8, the Poincare constant in \eqref{eq:spatialDiscontBound}, showing that \eqref{eq:spatialDiscontBound} is satisfied.}
\end{figure}

If we consider this problem as instance of a 1D nonlocal elastic bar where a portion of the bar has been excised, the region near the excised portion will have a reduced density of bonds, which will manifest as a locally reduced stiffness. Inspecting Figure~\ref{fig:horizon}, we see that the more material that is removed (i.e., the larger the value of $\varepsilon$) the steeper the slope of the displacement for the remaining material around $x=0.5$, which is consistent with a reduced stiffness in that region. 

Note from Theorem \ref{thm:kernelBound} we have if $\mu_i\in L^2(\ren\times\ren)$ and $\|\wt{\mu}_{1,\asym}\|_{L^2(\ren\times\ren)}<C_P^{-1}$, then
\[
    \|u_2-u_1\|_{L^2(\dom)}^2\le\frac{C_P}{1-C_P\|\wt{\mu}_{1,\asym}\|_{L^2(\ren\times\ren)}}
    \ls \|\wt{\mu}_2-\wt{\mu}_1\|_{L^2(\ren\times\ren)}\|u_2\|_{L^2(\dom\cup\bnd)}
    +
   \|K\|_{L^\infty(\dom)}\|f\|_{L^2(\dom)}\rs.
\]
If we let $\mu_1=\mu_{\varepsilon=0}=3\delta^{-3}$ (with no bonds removed), we satisfy the restriction as $\mu_{1,asym}=0$. Further, since $f=0$, we have the bound
\begin{equation} \label{eq:spatialDiscontBound}
    \|u_2-u_1\|_{L^2(\dom)}^2\le C_P \|u_2\|_{L^2(\dom\cup\bnd)}
    \|\wt{\mu}_2-\wt{\mu}_1\|_{L^2(\ren\times\ren)}.
\end{equation}

\subsection{Nonlinear Forcing}\label{sec:nonLinear}
We consider a nonlinear variation of \eqref{eq:nonlocalLaplacian} with right-hand side
\begin{equation}\label{eq:nonlinear1}
f_{\eta, \theta}(x,u)=2\frac{\eta\arctan{u}+\theta}{x^2+1},
\end{equation}
parameterized by $\eta,\theta>0$, and study the sensitivity of solutions to perturbations of these parameters. This is a slightly altered variation of the example from \cite[Section 5.1]{fossraduwright}, for which well-posedness and regularity of solution follow with simple alterations from the arguments presented in \cite{fossraduwright}. We choose the kernel $\mu(x,y)=3\delta^{-3}$, let $\delta = 0.2$, and set the collar condition $g(x)=0$.

\subsubsection{Sensitivity to perturbations in $\eta$}\label{sec:nonLinear1}
In this example we fix $\theta=1$ and consider perturbations in the solution by varying $\eta$. We denote $\eta_1=0$ as the unperturbed solution (i.e., $u_{\eta_1}(x)$ denotes the solution for $\eta=0$). From Theorem \ref{nonlinear3}, we know that a change in the nonlinearity induces a variation in the solution with a bound given by
\begin{equation}\label{eq:nonlinearbound1}
    \| u_2-u_1\|_{L^2(\Omega \cup \Gamma)} \leq  C\|f_2-f_1\|_{L^\infty(\Omega \cup \Gamma)},
\end{equation}
where $C=\displaystyle\frac{C_P}{1-M_1M_2C_P\|\gamma_{\text{asym}}\|_{L^\infty(\Omega)}-C_PL_2 }$. Since $\mu(x,y)$ is symmetric, the constant reduces to $C=\displaystyle\frac{C_P}{1-C_PL_1}$. The condition $C_P^{-1}>L_1=2\eta_1=0$ ($\eta_1=0$ removes the nonlinearity in the right-hand-side) is clearly satisfied and so $C=C_P$. In Figure~\ref{fig:nonlinear1} we plot solutions for various $\eta$, and in Table \ref{table:nonlinear1} we show numerically that \eqref{eq:nonlinearbound1} is satisfied for this example. 

\subsubsection{Sensitivity to perturbations in $\theta$}\label{sec:nonLinear2}
\begin{figure}
    \centering
    \begin{subfigure}[b]{0.5\textwidth}
      \centering
      \includegraphics[scale=0.12]{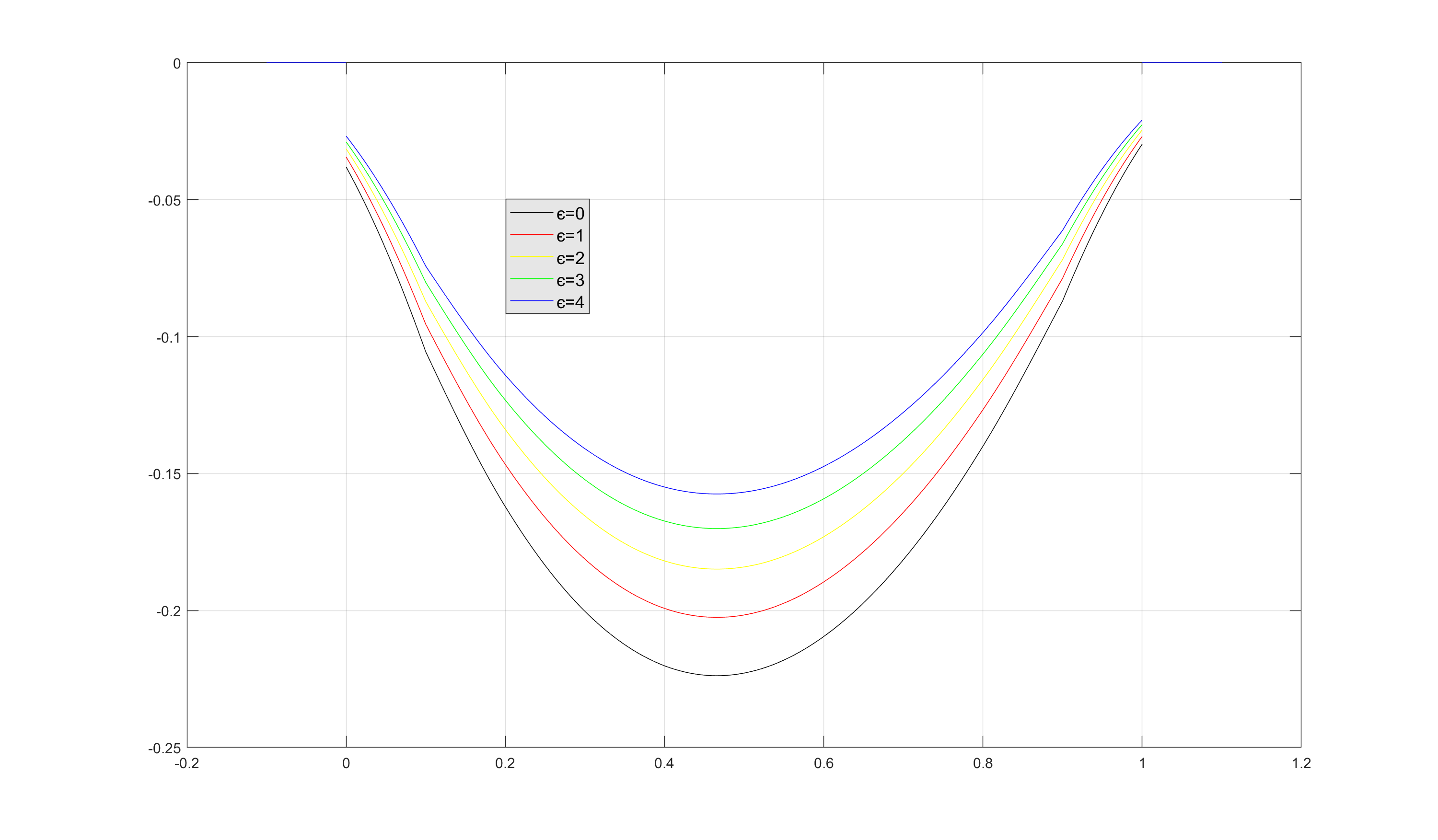}
      \caption{}
      \label{fig:nonlinear1}      
    \end{subfigure}
    \hfill
    \begin{subfigure}[b]{0.45\textwidth}
      \centering
      \includegraphics[scale=0.12]{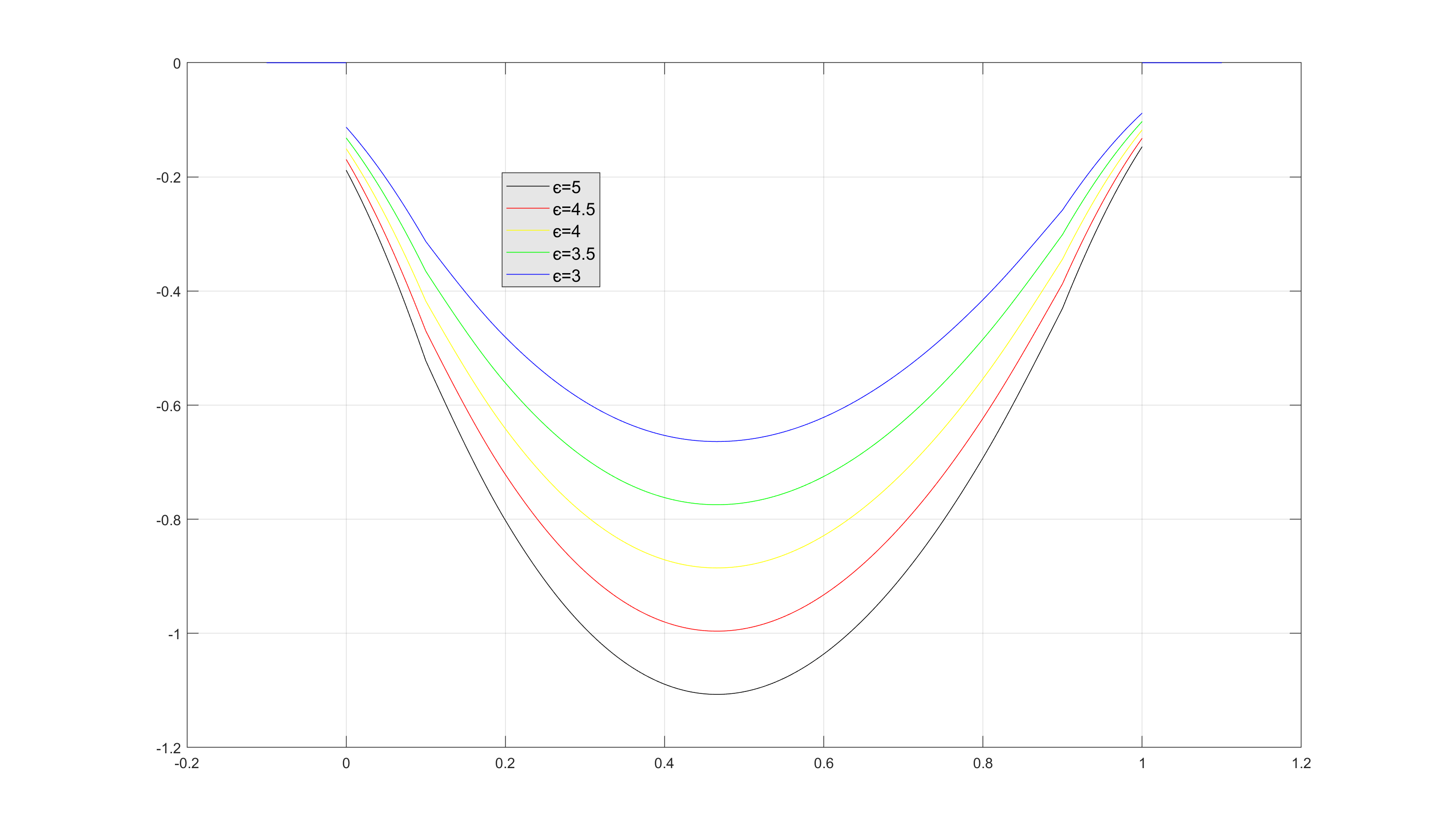}
      \caption{}
      \label{fig:nonlinear2}            
    \end{subfigure}
    \par \bigskip
    \begin{subtable}[b]{0.45\textwidth}
      \centering
      \begin{tabular}{|c|c|c|c|}
        \hline
         $\eta_2$	&$\|f_{\eta_2}-f_{\eta_1} \|_{L^\infty}$& $\|u_{\eta_2}-u_{\eta_1} \|_{L^2}$ &	 $\frac{\|u_{\eta_2}-u_{\eta_1} \|_{L^2}}{\|f_{\eta_2}-f_{\eta_1} \|_{L^\infty}}$ \\
\hline
1.0 &3.14159 &0.0159504 &0.00507719 \\
\hline
2.0 &6.28319 &0.0291401 &0.00463779 \\
\hline
3.0 &9.42478 &0.0402251 &0.00426802 \\
\hline
4.0 &12.5664 &0.0496702 &0.00395263 \\
        \hline
      \end{tabular}
    \caption{} 
    \label{table:nonlinear1}
    \end{subtable}
    \hfill
    \begin{subtable}[b]{0.45\textwidth}
      \centering
      \begin{tabular}{|c|c|c|c|}
        \hline
        $\theta$	&$\|f_{\theta_2}-f_{\theta_1} \|_{L^\infty}$& $\|u_{\theta_2}-u_{\theta_1} \|_{L^2}$ &	 $\frac{\|u_{\theta_2}-u_{\theta_1} \|_{L^2}}{\|f_{\theta_2}-f_{\theta_1} \|_{L^\infty}}$ \\
\hline
4.5&0.5 &0.0829955  &0.165991  \\
\hline
4.0& 1&0.165959 &0.165959  \\
\hline
3.5&1.5  &0.248893 &0.165928  \\
\hline
3.0&2 &0.331797 &0.165899   \\
        \hline
      \end{tabular}
    \caption{} 
    \label{table:nonlinear2}
    \end{subtable}
    \caption{Numerical results for the examples of \S\ref{sec:nonLinear1} and \S\ref{sec:nonLinear2}. (a): Solutions $u_\eta(x)$ over $(-\delta,1+\delta)$ for various $\eta$ using the nonlinear forcing function \eqref{eq:nonlinear1}. (b): Solutions $u_\theta(x)$ over $(-\delta,1+\delta)$ for various $\theta$ using the nonlinear forcing function \eqref{eq:nonlinear1}.  (c): Numerical evaluation of the terms in \eqref{eq:nonlinearbound1} for various $\eta$. (d): Numerical evaluation of the terms in \eqref{eq:nonlinearbound1} for various $\theta$.}
    \label{fig:nonlinear}        
\end{figure}

In this example we fix $\eta=1/9$, and consider perturbations in the solution by varying $\theta$. We denote $\theta_1=5$ as the unperturbed solution (i.e., $u_{\theta_1}(x)$ denotes the solution for $\theta=5$). Referring again to the bound \eqref{eq:nonlinearbound1}, since $\mu(x,y)$ is symmetric, the constant in this example reduces to $C=\displaystyle\frac{C_P}{1-C_PL_1}$. The condition $C_P^{-1}>L_1=2\eta_1=2/9$, and thus $C = 3/2$. In Figure~\ref{fig:nonlinear2} we plot solutions for various $\theta$, and in Table \ref{table:nonlinear2} we show numerically that \eqref{eq:nonlinearbound1} is satisfied for this example. 

\subsubsection{An exponential kernel} \label{sec:expKernel}
\begin{figure}
    \centering
    \begin{subfigure}[b]{0.45\textwidth}
      \includegraphics[scale=0.12]{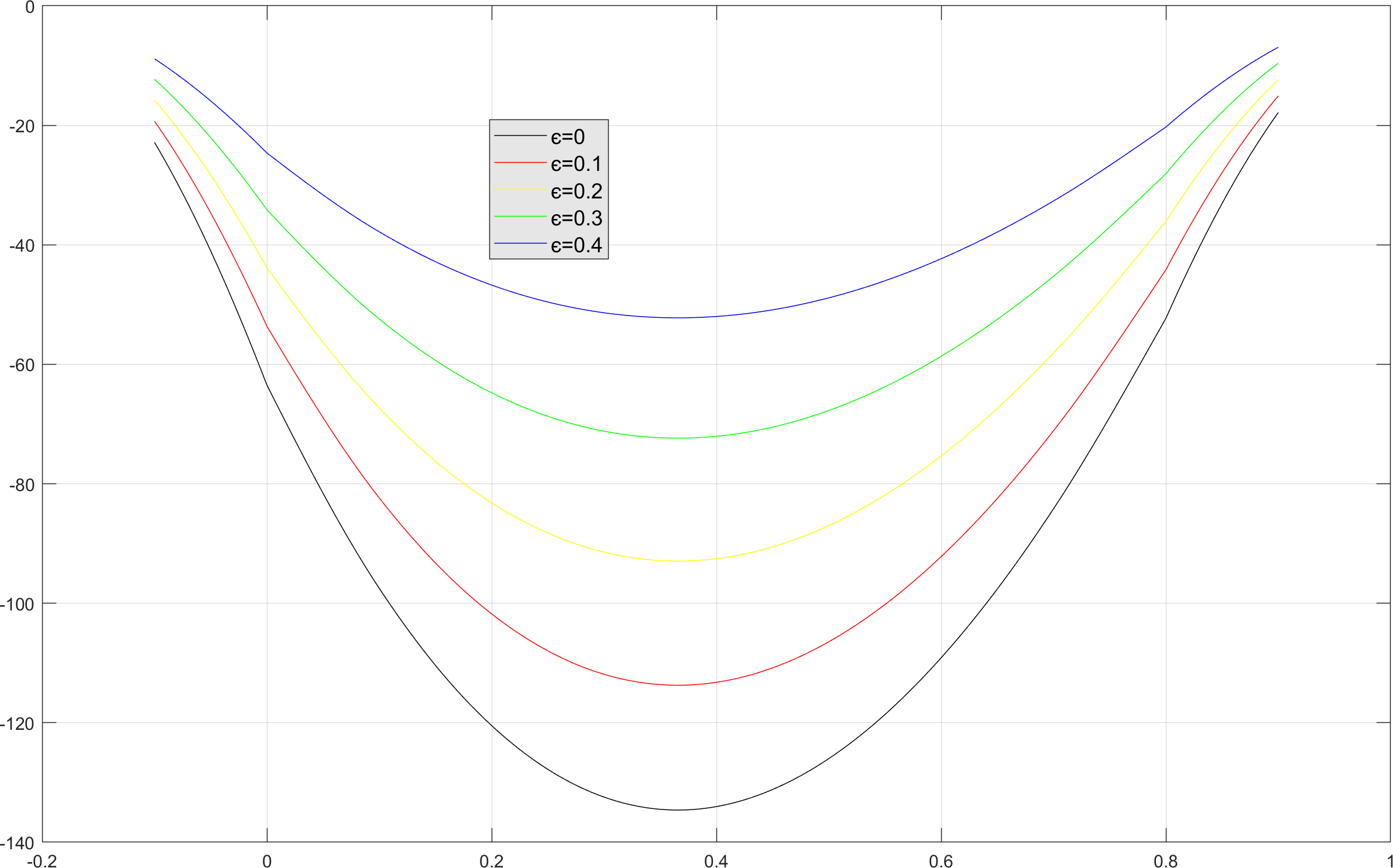}
      \caption{}
      \label{fig:exp_nonlinear1}      
    \end{subfigure}
    \hfill
    \centering
    \begin{subfigure}[b]{0.45\textwidth}
      \includegraphics[scale=0.12]{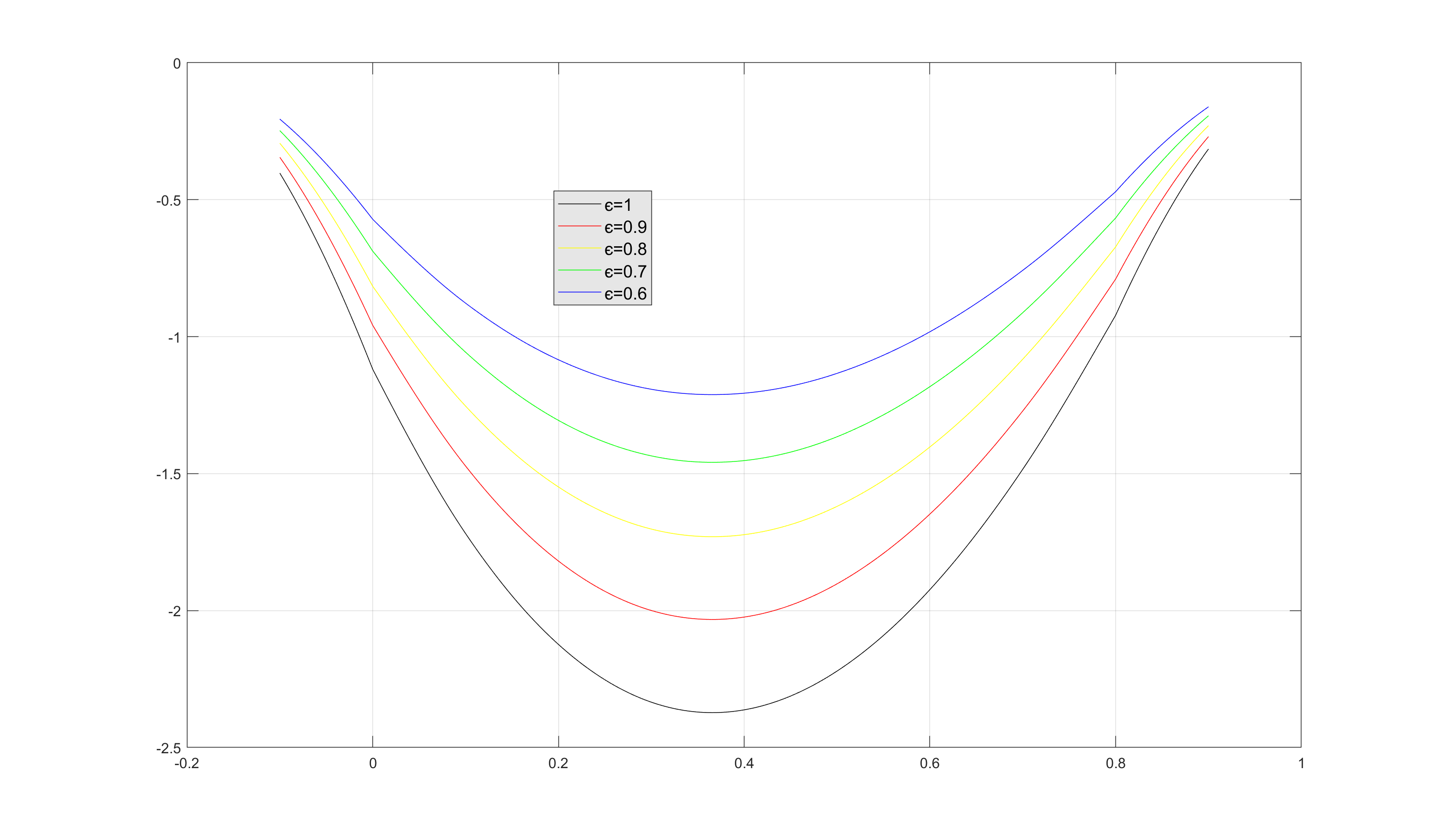}
      \caption{}
      \label{fig:exp_nonlinear2}            
    \end{subfigure}
    \par \bigskip
    \begin{subtable}[b]{0.45\textwidth}
      \begin{tabular}{|c|c|c|c|}
        \hline
         $\eta_2$	&$\|f_{\eta_2}-f_{\eta_1} \|_{L^\infty}$& $\|u_{\eta_2}-u_{\eta_1} \|_{L^2}$ &	 $\frac{\|u_{\eta_2}-u_{\eta_1} \|_{L^2}}{\|f_{\eta_2}-f_{\eta_1} \|_{L^\infty}}$ \\
\hline
0.1 &0.314159 &15.6436 & 49.7952 \\
\hline
0.2 & 0.628319&31.2032 & 49.6615 \\
\hline
0.3 &0.942478 &46.6088 & 49.4534 \\
\hline
0.4 & 1.25664&61.69 &49.0913  \\
        \hline
      \end{tabular}
    \caption{} 
    \label{table:exp_nonlinear1}
    \end{subtable}
    \hfill
    \begin{subtable}[b]{0.45\textwidth}
      \begin{tabular}{|c|c|c|c|}
        \hline
        $\theta$	&$\|f_{\theta_2}-f_{\theta_1} \|_{L^\infty}$& $\|u_{\theta_2}-u_{\theta_2} \|_{L^2}$ &	 $\frac{\|u_{\theta_2}-u_{\theta_1} \|_{L^2}}{\|f_{\theta_2}-f_{\theta_1} \|_{L^\infty}}$ \\
\hline
0.9&0.00505439  &0.254709 &50.3935  \\
\hline
0.8& 0.00954429 &0.48097 &50.3935  \\
\hline
0.7&0.0135794  &0.684312 &50.3935  \\
\hline
0.6&0.0172518  &0.869377 &50.3935  \\
        \hline
      \end{tabular}
    \caption{} 
    \label{table:exp_nonlinear2}
    \end{subtable}
    \caption{Numerical results for the example of \S\ref{sec:expKernel}. This example uses the exponential kernel of \eqref{gauss}; compare against solutions using a constant kernel in Figure~\ref{fig:nonlinear}. (a): Solutions $u_\eta(x)$ over $(-\delta,1+\delta)$ for various $\eta$ using the nonlinear forcing function \eqref{eq:nonlinear1} with the exponential kernel. (b): Solutions $u_\theta(x)$ over $(-\delta,1+\delta)$ for various $\theta$ using the nonlinear forcing function \eqref{eq:nonlinear1} with the exponential kernel. (c): Numerical evaluation of the terms in \eqref{eq:nonlinearbound1} for various $\eta$. (d): Numerical evaluation of the terms in \eqref{eq:nonlinearbound1} for various $\theta$.}
    \label{fig:expKernel}    
\end{figure}

The magnitude of the constant $C$ in bounds of the form \eqref{eq:nonlinearbound1} depends on the conditioning of the operator, which depends strongly on the choice of kernel. In this example we repeat the analysis of \S\ref{sec:nonLinear1} and \S\ref{sec:nonLinear2} with the kernel 
\begin{equation}\label{gauss}
\mu(x,y)=
\begin{cases}
c_{\delta}e^{-(x-y)^2}, & |x-y|<\delta\\
0, & |x-y|\geq \delta,
\end{cases}
\end{equation} where $c_{\delta}$ is chosen such that $\|\mu\|_{L^1(\re)}=1$. Results are shown in Figure~\ref{fig:expKernel}. Specifically, observe the magnitude of the solution in Figures~\ref{fig:exp_nonlinear1} and~\ref{fig:exp_nonlinear2} is substantially greater than in Figures~\ref{fig:nonlinear1} and~\ref{fig:nonlinear1}, which is consistent with the rightmost columns in Tables~\ref{table:exp_nonlinear1} and~\ref{table:exp_nonlinear2} being substantially larger than the rightmost columns of Tables~\ref{table:nonlinear1} and~\ref{table:nonlinear2}. 

\section{Conclusions and Future Work}

The results proven show explicit dependence of solutions to linear and nonlinear nonlocal systems with respect to forcing terms (including nonlinear Lipschitz forcing), Dirichlet boundary conditions, and different choices for kernels. In the case of heterogeneous kernels, additional restrictions are needed in order to accommodate the explicit dependence on space variables. 

The numerical studies validate the theoretical bounds, which are based on upper bounds for the Poincar\'e constant. The simulations, however, seem to suggest that the theoretical bounds obtained for the sensitivity results are not optimal, hinting at an open research direction. 

Generalizing the results to the vector valued framework, as given by the state-based peridynamics formulation would be important for a variety of applications which use this theory. Additionally, real-world applications may require different types of nonlinearities, so eliminating the Lipschitz restriction would provide a significant advance for the stability of these problems. Future work will consider Neumann (or flux-type), as well as mixed-type (possibly nonlinear) boundary conditions are an expected future step in understanding the effect that data imposed on collar (however small) may have on solutions. Finally, we are exploring stability results for higher-order systems, in particular, involving a nonlocal  biharmonic operator.

\section{Acknowledgements}
MLP was supported by the U.S. Department of Energy, Office of Science, Office of Advanced Scientific Computing Research under the Physics-Informed Learning Machines for Multiscale and Multiphysics Problems (PhILMs) project.Sandia National Laboratories is a multimission laboratory managed and operated by National Technology \& Engineering Solutions of Sandia, LLC, a wholly owned subsidiary of Honeywell International Inc., for the U.S. Department of Energy’s National Nuclear Security Administration under contract DE-NA0003525. This paper describes objective technical results and analysis. Any subjective views or opinions that might be expressed in the paper do not necessarily represent the views of the U.S. Department of Energy or the United States Government.

\bibliographystyle{siam}
\bibliography{NumCont}

\end{document}